%% file: sc17askit_acm.tex
\documentclass[sigconf,screen,natbib]{acmart}

\input preamble.tex

\input macros.tex

\input notation.tex

\begin{document}

\title{Geometry-Oblivious FMM for Compressing  Dense SPD Matrices}

\author{Chenhan~D. Yu}
\affiliation{%
  \institution{{\wu}Department of Computer Science}
  \institution{{\wg}Institute for Computational Engineering and Sciences}
  \streetaddress{The University of Texas at Austin, Austin, TX}}
\email{chenhan@cs.utexas.edu}

\author{James Levitt}
\affiliation{%
  \institution{Institute for Computational Engineering and Sciences}
  \streetaddress{The University of Texas at Austin, Austin, TX}}
\email{jlevitt@ices.utexas.edu}

\author{Severin Reiz}
\affiliation{%
  \institution{Institute for Computational Engineering and Sciences}
  \streetaddress{The University of Texas at Austin, Austin, TX}}
\email{s.reiz@tum.de}

\author{George Biros}
\affiliation{%
  \institution{Institute for Computational Engineering and Sciences}
  \streetaddress{The University of Texas at Austin, Austin, TX}}
\email{gbiros@acm.org}

\begin{abstract}
\input abstract.tex
\end{abstract}

\maketitle

\section{Introduction} \label{s:intro} \input intro.tex  

\section{Methods} \label{s:methods} \input methods.tex   %
  
  \subsection{Geometry-oblivious techniques} \label{s:oblivious} \input geomOblivious.tex   

  \subsection{Algebraic Fast Multipole Method} \label{s:algfmm} \input algfmm.tex   

  \subsection{Shared memory parallelism} \label{s:par} \input par.tex   

\section{Experimental Setup} \label{s:setup} \input setup.tex   

\section{Empirical Results} \label{s:results} \input results.tex 

\section{Conclusions} \label{s:conclusion} \input conclusion.tex   

\newpage

\bibliographystyle{ACM-Reference-Format}
\bibliography{gb,refs}

\newpage
\section*{Appendix} \input appendix.tex 

\end{document}

%% file: preamble.tex
\usepackage{graphics,epsfig,graphicx,float,subfigure,color}
\usepackage{algorithmic}
\usepackage[margin=0pt,labelsep=space,font=small,labelfont=bf,textfont=it]{caption}
\usepackage{atbegshi}

\usepackage[section]{algorithm}
\usepackage{lineno}
\usepackage{verbatimbox}

\usepackage{etex}

\usepackage{amsmath,amssymb,amsbsy,amsfonts,latexsym}
\usepackage{multicol,multirow}
\usepackage{hhline}
\usepackage{paralist}
\usepackage{array}
\usepackage{marvosym}
\usepackage{url}
\usepackage{boxedminipage}
\usepackage{wrapfig}
\usepackage{multirow}
\usepackage{threeparttable}
\usepackage{pdfpages}
\usepackage{xspace}
\usepackage{colortbl}

\usepackage{tabularx,booktabs}

\usepackage{xcolor}

\usepackage{tikz} 
\usetikzlibrary{arrows} 
\usetikzlibrary{patterns}  
\usepackage{pgfplots} 
\pgfplotsset{compat=1.3}

\let\labelindent\relax
\usepackage{enumitem}

\usepackage{siunitx}
\usepackage{booktabs}
\usepackage{calligra}

\DeclareMathAlphabet{\mathcalligra}{T1}{calligra}{m}{n}

\definecolor{light-gray}{gray}{0.8}
\newcommand{\accnum}{\num[scientific-notation=true, round-mode=places, round-precision=0]}
\newcolumntype{R}{>{\columncolor{light-gray}}r}
\newcolumntype{L}{>{\columncolor{light-gray}}l}
\newcolumntype{C}{>{\columncolor{light-gray}}c}

%% file: macros.tex
\def\gap{\hspace*{.2in}}

\newcommand{\figref}[1]{Figure~\ref{#1}}
\newcommand{\tabref}[1]{Table~\ref{#1}}
\newcommand{\secref}[1]{\S\ref{#1}}
\newcommand{\algref}[1]{Algorithm~\ref{#1}}

\newcommand{\MA}[1]{{\mathcal #1}}
\DeclareMathOperator{\bigO}{\mathcal{O}}

\newcommand{\zapspace}{\topsep=0pt\partopsep=0pt\itemsep=0pt\parskip=0pt}

\definecolor{light-gray}{gray}{0.80}

\newcommand{\ipoint}[1]{{\color{darkgray}\textit{\textbf{#1}}}}

\newcounter{magicrownumbers}
\newcommand\rownumber{\refstepcounter{magicrownumbers}\arabic{magicrownumbers}}

%% file: notation.tex
\newcommand{\Ker}{\mathcal{K}}

\newcommand{\sk}[1]{{\widetilde{#1}}} 
\newcommand{\K}{{\MA{K}}} 
\newcommand{\lc}{{\tt l}}  
\newcommand{\rc}{{\tt r}}  
\newcommand{\ns}{{s}} 

\newcommand{\idtol}{{\tau}}

\newcommand{\hmatrix}{\ensuremath{\mathcal{H}\!}-matrix}
\newcommand{\gofmm}{{\texttt{GOFMM}}}

\newcommand{\superscript}[1]{\ensuremath{^{\textrm{#1}}}}

\def\wu{\superscript{*}}
\def\wg{\superscript{\dag}}

%% file: abstract.tex
We present \gofmm{} (geometry-oblivious FMM), a novel method that creates a
hierarchical low-rank approximation, or \emph{``compression,''} of an arbitrary
dense symmetric positive definite (SPD) matrix. For many applications, \gofmm{}
enables an approximate matrix-vector multiplication in $N \log N$ or even $N$
time, where $N$ is the matrix size. Compression requires $N \log N$ storage and
work.  In general, our scheme belongs to the family of hierarchical matrix
approximation methods. In particular, it generalizes the fast multipole method
(FMM) to a purely algebraic setting by only requiring the ability to sample
matrix entries. Neither geometric information (i.e., point coordinates) nor
knowledge of how the matrix entries have been generated is required,  thus the
term \emph{``geometry-oblivious.''} Also, we introduce a shared-memory parallel scheme for hierarchical matrix computations that reduces synchronization barriers. We present results on the Intel Knights Landing and Haswell architectures, and on the NVIDIA Pascal architecture for a variety of matrices.

%% file: intro.tex
\begin{figure}[!t]
  \centering
  \includegraphics[scale=.27]{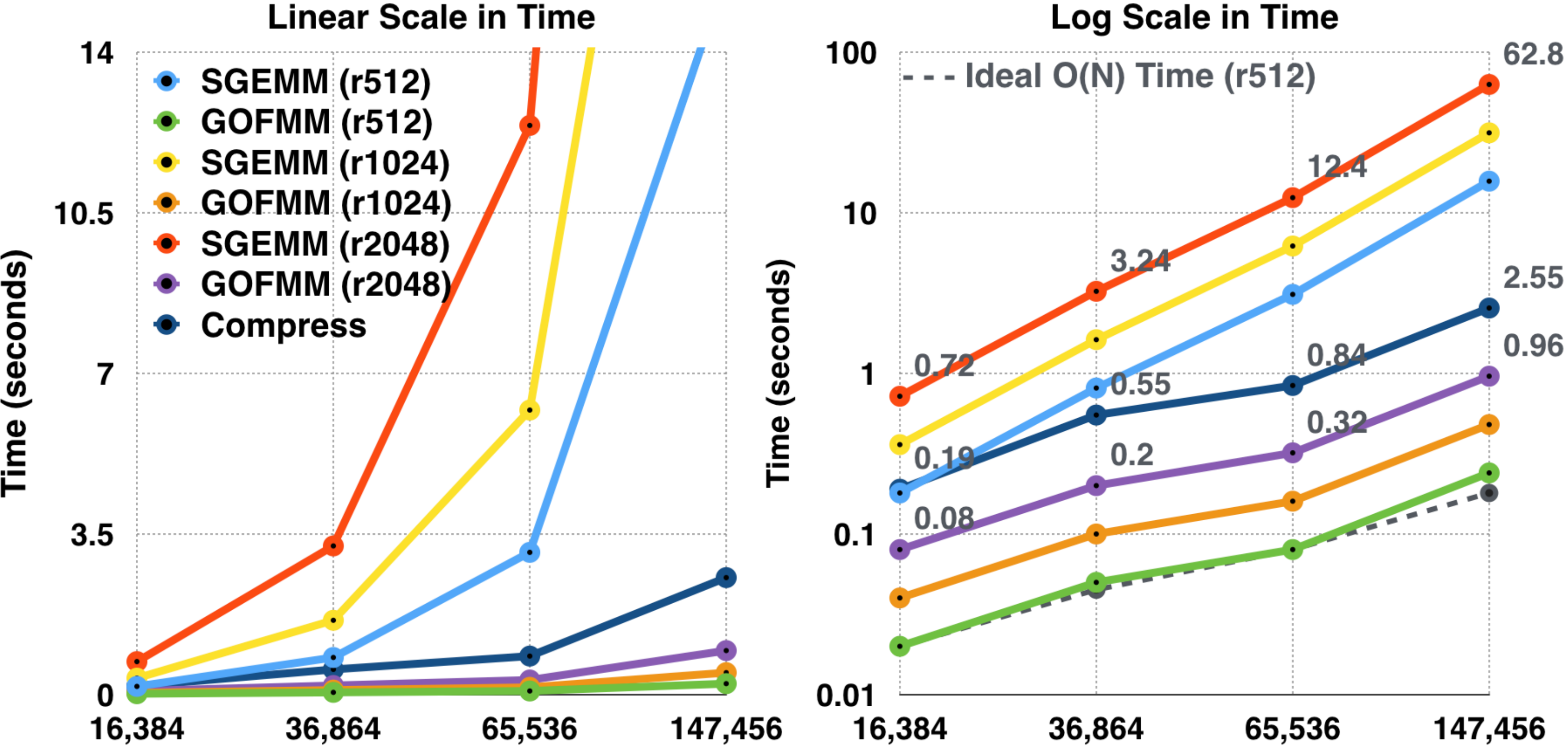}
  \caption{
					Comparison of \ipoint{runtime in seconds (y-axis)} versus \ipoint{problem
					size $N$ (x-axis)} to multiply test matrix K02 (see \secref{s:setup}) of size $N \times N$
	with a matrix of size $N \times r$, for $r = 512, 1024, 2048$.
	Results are plotted against a linear scale (left) and a logarithmic scale (right).
	The top three curves demonstrate $\MA{O}(N^2)$ scaling of Intel MKL \texttt{SGEMM} for
	each value of $r$.
	The middle curve shows the time for \gofmm{} to compress K02,
	which scales as $\MA{O}(N\log{N})$ in these cases.
	The bottom three curves show the $\MA{O}(N)$ scaling of the time for \gofmm{}
	to evaluate the matrix product for each value of $r$
	after compression is already completed.
	The \gofmm{} results reach accuracies of \num{1e-2} to \num{4e-4} in single precision.
	In these experiments, the crossover (including compression time) is $N=\num{16384}$,
	and for $N=\num{147456}$, we observe an 18$\times$ speedup over \texttt{SGEMM}.
	}
  \label{fig:quadratic}
\end{figure}

We present \gofmm{}, a novel  algorithm for the approximation of dense symmetric positive definite (SPD) matrices. \gofmm{} can be used for compressing a dense matrix and accelerating matrix-vector multiplication operations. As an example, in~\figref{fig:quadratic} we report timings for an SGEMM operation using an optimized dense matrix library and compare with the \gofmm{}-compressed version.

Let $K \in \mathbb{R}^{N\times N}$ be a dense SPD matrix, with $K=K^T$ and $x^T
K x > 0,\ \forall x\in \mathbb{R}^N,\ x \neq 0$. Since $K$ is dense it requires $\bigO(N^2)$ storage  and $\bigO(N^2)$ work for a matrix-vector multiplication (hereby \ipoint{``matvec''}). Using $\bigO(N\log N)$ memory and work, we construct an approximation $\tilde{K}$ such that $\|\tilde{K}-K\| \leq \epsilon \|K\|$, where $\epsilon$ is a user-defined error tolerance. Assuming the evaluation of a single matrix entry $K_{ij}$ requires $\bigO(1)$ work, a matvec  with $\tilde{K}$  requires $\bigO(N \log N)$ or $\bigO(N)$ work depending on the properties of $K$ and the \gofmm{} variant. Our scheme belongs to the class of hierarchical matrix approximation methods.

{\bf Problem statement:} given any SPD matrix $K$, our task is to construct a hierarchically low-rank matrix $\tilde{K}$ so that $\|K-\tilde{K}\|/\|K\|$ is small.  The only required input to our algorithm is a routine that returns $K_{IJ}$, for arbitrary row and column index sets $I$ and $J$. The constant in the complexity estimate depends on the structure of the underlying matrix. Let us remark and emphasize that our scheme \emph{cannot guarantee both accuracy and work complexity} simultaneously since an arbitrary SPD matrix may not admit a good hierarchical matrix  approximation (see ~\secref{s:methods}). 

We say that a matrix $\tilde{K}$ has a \ipoint{hierarchically low-rank structure}, i.e., $\tilde{K}$ is an \hmatrix{}~\cite{hackbusch15,bebendorf08}, if 
\begin{equation}\label{e:hmatrix}
\tilde{K}=D+S+UV,
\end{equation}
where $D$ is \ipoint{block-diagonal} with \ipoint{every block being an \hmatrix{}}, $U$ and $V$  are
\ipoint{low rank}, and $S$ is \ipoint{sparse}.  At the base case of this
recursive definition, the blocks of $D$ are small dense matrices. An \hmatrix{}
matvec requires $\bigO(N \log N)$ work the constant depending on the rank of
$U$ and $V$. Depending on the construction algorithm, this complexity can go down to $\bigO(N)$. Although such matrices are rare in real-world  applications, it is quite common to find matrices that  can be approximated \emph{arbitrarily} well by an \hmatrix{}.

One important observation is that \emph{this hierarchical low-rank structure  is not invariant to row and column permutations}. Therefore any algorithm for constructing $\tilde{K}$ must first appropriately permute $K$ before constructing the matrices $U, V, D$, and $S$. Existing algorithms  rely on the matrix entries $K_{ij}$ being ``interactions'' (pairwise functions) between \ipoint{points} $\{x_i\}_{i=1}^N$ in $\mathbb{R}^d$ and permute $K$ either by clustering the points (typically using some tree data-structure) or by using graph partitioning techniques (if $K$ is sparse). \gofmm{} does not require such geometric information. 

\paragraph{\textbf{Background and significance}}
Dense SPD matrices appear in scientific computing, statistical inference, and data analytics. They appear in Cholesky and LU factorization~\cite{grasedyck-kriemann-leborne08}, in Schur complement matrices for saddle point problems~\cite{benzi-golub-liesen05}, in Hessian operators in optimization~\cite{biegler-ghattas01}, in kernel methods for statistical learning~\cite{hofmann-scholkopf-smola08,gray-moore01}, and in N-body methods and integral equations~\cite{greengard94,hackbusch15}. In many applications, the entries of the input matrix $K$ are given by $K_{ij} = \Ker(x_i, x_j):\mathbb{R}^d\times\mathbb{R}^d\rightarrow \mathbb{R}$, where $\Ker$ is a \ipoint{kernel function}. Examples of kernel functions are radial basis functions, Green's functions, and angle similarity functions. For such \emph{kernel matrices}, the input is not a matrix, but only the points $\{x_i\}_{i=1}^N$. The  points are used to appropriately permute the matrix using spatial data structures. Furthermore, the construction of the sparse correction $S$ uses  nearest-neighbor structure of the input points. The low-rank matrices $U,V$ can be either analytically computed using expansions of the kernel function, or semi-algebraically computed using fictitious points (or equivalent points), or using algebraic sampling-based methods that use geometric information. In a nutshell, geometric information is used in all aspects of an \hmatrix{} method.

In many cases however, such points and kernel functions are not available. For example, in dense graphs in data analysis (e.g., social networks, protein interactions). Related matrices include graph Laplacian operators and their inverses. Additional examples include frontal matrices and Schur complements in factorization of sparse matrices; Hessian operators in optimization; and kernel methods in machine learning without points (e.g., word sequences and diffusion on graphs~\cite{cancedda-e03,kondor-lafferty02}).

\paragraph{\textbf{Contributions}}
\gofmm{} is inspired by the rich literature of algorithms for matrix sketching, hierarchical matrices, and fast multipole methods.  Its unique feature is that by using only matrix evaluations it generalizes FMM ideas to compressing arbitrary SPD matrices.  In more detail, our contributions are summarized below.
\begin{itemize}[leftmargin=*]\zapspace
  \item A result from reproducing kernel Hilbert space theory is that any SPD
    matrix corresponds to a Gram matrix of vectors in some, unknown Gram (or
    feature) space~\cite{hofmann-scholkopf-smola08}. Based on this result, the matrix entries are inner products, which we use to define distances. These distances allow us to design an efficient, purely algebraic FMM method. 
  \item The key algorithmic components of \gofmm{} (and other hierarchical
    matrix and FMM codes) are tree traversals. 
    We test parallel level-by-level traversals, \emph{out-of-order} traversals
    using \texttt{OpenMP}'s advanced task scheduling and an in-house tree-task scheduler. 
    We found that scheduling significantly improves the performance when compared to 
    level-by-level tree traversals. We also use this scheduling to support heterogeneous architectures. 
  \item We conduct extensive experiments to demonstrate the feasibility of the
    proposed approach. We test our code on 22 different matrices related to
    machine learning, stencil PDEs, spectral PDEs, inverse problems, and graph
    Laplacian operators. We perform numerical experiments on Intel Haswell
    and KNL, Qualcomm ARM, and NVIDIA Pascal architectures. 
    Finally, we compare with three state-of-the-art codes: \texttt{HODLR},
    \texttt{STRUMPACK}, and \texttt{ASKIT}. 
\end{itemize}
\gofmm{} also has several additional capabilities. If points and kernel
functions (or Green's) function are available, they can be utilized in a
similar way to the algebraic FMM code \texttt{ASKIT} described
in~\cite{march-xiao-yu-biros15,march-xiao-biros-fmm-e15}. \gofmm{} currently
supports three different measures of distance: geometric point-based (if
available), Gram-space $\ell^2$ distance, and Gram-space angle distance.  \gofmm{} has support for matvecs with \ipoint{multiple vectors}, which is useful for Monte-Carlo sampling, optimization, and block Krylov methods.

\paragraph{\textbf{Limitations}}
\gofmm{} is restricted to SPD matrices. (However, if we are given points, the method becomes similar to existing methods).
\gofmm{} guarantees symmetry of $\tilde{K}$, but if $\|K-\tilde{K}\|/\|K\|$ is large, positive definiteness may be compromised. 
To reiterate, \gofmm{} cannot simultaneously guarantee both accuracy and work complexity.
This initial implementation of \gofmm{} supports shared-memory parallelism and accelerators, but not distributed memory architectures.
The current version of \gofmm{} also has several parameters that require manual tuning.
Often, the main goal of building \hmatrix{} approximations is to construct a factorization of $K$, a topic we do not discuss in this paper.  Our method requires the ability to evaluate kernel entries and the complexity estimates require that these entries can be computed in $\bigO(1)$ time. If $K$ is only available through matrix-free interfaces, these assumptions may not be satisfied. Other algorithms, like \texttt{STRUMPACK}, have inherent support for such matrix-free compression.

\paragraph{\textbf{Related work.}} The literature on hierarchical matrix methods and fast multipole methods is vast. Our discussion is brief and limited to the most related work.
                              
\ipoint{Low-rank approximations.} The most popular approach for compressing arbitrary matrices is a global low-rank approximation using randomized linear algebra. In~\eqref{e:hmatrix}, this is equivalent to setting $D$ and $S$ to zero and constructing only $U$ and $V$. Examples include the CUR~\cite{mahoney-drineas09} factorization, the Nystrom approximation~\cite{williams-seeger01}, the adaptive cross approximation~\cite{bebendorf-rjasanow03}, and randomized rank-revealing factorizations~\cite{martinsson-rokhlin-tygert10,halko-martinsson-tropp11}. These techniques can also be used for \hmatrix{} approximations when $D$ is not zero. Instead of applying them to $K$, we can apply them to the off-diagonal blocks of $K$. FMM-specific techniques that are
a mix between analytic and algebraic methods include
kernel-independent methods~\cite{martinsson-rokhlin07,ying-biros-zorin-03} and the
black-box FMM~\cite{fong-darve09}. Constructing both $U$ and $V$ accurately and with optimal complexity is hard. The most robust algorithms require $\bigO(N^2)$ complexity or higher (randomized methods and leverage-score sampling) since they require one to ``touch'' all the entries of the matrix (or block) to be approximated.

\ipoint{Permuting the matrix.} When $K$ is sparse, the method of choice uses graph-partitioning. This doesn't scale to dense matrices because practical graph partitioning  algorithms scale at least linearly with the number of edges and thus the construction cost would be at least $\bigO(N^2)$~\cite{agullo-darve-e16,karypis-kumar98}. 

\ipoint{\hmatrix{} methods and software.} Treecodes and fast multipole methods originally were developed for N-body problems and integral equations. Algebraic variants led the way to the abstraction of \hmatrix{} methods and the application to the factorization of sparse systems arising from the discretization of elliptic PDEs~\cite{hackbusch15,bebendorf08,ambikasaran-13,greengard-gueyffier-martinsson-rokhlin09,ho-greengard12,xia-e10}. 

Let us briefly summarize the \hmatrix{} classification. Recall the decomposition $K=D+UV+S$, \eqref{e:hmatrix}. If $S$ is zero the approximation is called a hierarchically off-diagonal low rank (HODLR) scheme.  In addition to $S$ being zero, if the \hmatrix{} decomposition of $D$ is used to construct $U$, $V$ we have a hierarchically semi-separable (HSS) scheme. If $S$ is not zero we have a generic \hmatrix{}; but if the $U,V$ terms are constructed in a nested way then we have an $\mathcal{H}^2$-matrix or an FMM depending on more technical details.  HSS and HODLR matrices lead to very efficient approximation algorithms for $K^{-1}$. However, $\MA{H}^2$ and FMM compression schemes better control the maximum rank of the $U$ and $V$ matrices  than HODLR and HSS schemes. For the latter, the rank of $U$ and $V$ can grow with $N$~\cite{chandrasekaran-gu-e10} and the complexity bounds are no longer valid. Recently, here have been algorithms to effectively compress FMM and $\MA{H}^2$-matrices~\cite{coulier-pouransari-darve16,yokota-ibeid-keyes16}. 
\begin{table}
\begin{tabular}{|>{\columncolor[gray]{0.8}}l|c|c|c|c|} 
\hline
\rowcolor[gray]{0.8}
\textbf{METHOD}  & \textbf{MATRIX} &  \textbf{LOW-RANK} & \textbf{PERM} & $S$ \\
\hline
\texttt{FMM}~\cite{cheng-greengard-rokhlin-99}
      &  $\Ker(x_i,x_j)$  &  EXP & OCTREE & Y \\
\texttt{KIFMM}~\cite{ying-biros-zorin-03}
      &  $\Ker(x_i,x_j)$  &  EQU & OCTREE & Y \\
\texttt{BBFMM}~\cite{fong-darve09}
      &  $\Ker(x_i,x_j)$  &  EQU & OCTREE & Y \\
\texttt{HODLR}~\cite{ambikasaran-darve13}
      &  $K_{ij}$       &  ALG & NONE & N \\
\texttt{STRUMPACK}~\cite{rouet-li-e16}
      &  $K_{ij}$       &  ALG & NONE & N \\            
\texttt{ASKIT}~\cite{march-xiao-yu-biros-sisc16}
&  $\Ker(x_i,x_j)$  & ALG  & TREE & Y \\
\texttt{MLPACK}~\cite{mlpack13}
& $\Ker(x_i,x_j)$ &  EQU  & TREE & Y \\
\textbf{\gofmm} &  $K_{ij}$     & ALG & TREE & Y \\
\hline
\end{tabular}
\caption{We summarize the main features of different \hmatrix{} methods/codes
for dense matrices. ``\textbf{MATRIX}'' indicates whether the method requires a
kernel function and points---indicated by $\Ker(x_i,x_j)$---or it just requires kernel entries---indicated by $K_{ij}$.  ``\textbf{LOW-RANK}'' indicates the method used for the off-diagonal low-rank approximations: ``EXP'' indicates kernel function-dependent analytic expansions; ``EQU'' indicates the use  of equivalent points (restricted to low $d$ problems); ``ALG'' indicates an algebraic method. ``\textbf{PERM}'' indicates the permutation scheme used for dense matrices: ``OCTREE'' indicates that the scheme doesn't generalize to high dimensions; ``NONE'' indicates that the input lexicographic order is used; and ``TREE'' indicates geometric partitioning that scales to high dimensions.  $S$ indicates whether a sparse correction (FMM or $\MA{H^2}$) is supported. In~\secref{s:results}, we present comparisons with \texttt{ASKIT}, STRUMPACK, and HODLR.}\label{t:codes}
\end{table}
One of the most scalable methods is \href{http://portal.nersc.gov/project/sparse/strumpack/}{STRUMPACK}~\cite{ghysels-li-e16,rouet-li-e16,martinsson16}, which constructs an HSS approximation of a square matrix (not necessarily SPD) and then uses it to construct an approximate factorization. For dense matrices STRUMPACK uses the lexicographic ordering. If no fast matrix-vector multiplication is available, STRUMPACK requires $\bigO(N^2)$ work for compressing a dense SPD matrix, and $\bigO(N)$ work for the matvec.

          

%% file: methods.tex

Given $K\in\mathbb{R}^{N\times N}$, \gofmm{} aims to construct an \hmatrix{}
$\sk{K}$ in the form of \eqref{e:hmatrix} such that we can approximate 
\begin{equation}
  u=Kw\approx \sk{K}w,\ \mbox{~for~}\ w\in\mathbb{R}^{N}.
  \label{e:appromatvec}
\end{equation}
When points $\{x_i\}^{N}_{i=1}$ are available such that $K_{ij} = \Ker(x_i,x_j)$,
the recursive partitioning on $D$ and the low-rank structure  $UV$ use \emph{distances} between $x_i$ and $x_j$.
Existing FMM methods approximate $K_{ij}$ when $x_i$ and $x_j$ are sufficiently \emph{far} from each other.
Otherwise, $K_{ij}$ is not approximated and it is placed either in $D$ or in $S$. 
We call this distance-based criterion \ipoint{near-far pruning}. 

To define such a pruning scheme without $\{x_i\}^{N}_{i=1}$,
we need a notion of distance between two matrix indices $i$ and $j$. We define such a distance in the next section.
With it, we can permute $K$ and define neighbors for each index $i$.
In \secref{s:algfmm}, we describe a task-based algebraic FMM that only relies on the distance we define.
Finally in \secref{s:par}, we discuss task parallelism and scheduling.


%% file: geomOblivious.tex
In this section, we introduce the machinery for using \gofmm{} in a
geometry-oblivious manner. Throughout the following discussion, we refer to a
set of indices $\mathcal{I} = \{1,\dots,N\}$, where index $i$ corresponds to
the $i$th row (or column) of the matrix $K$ in the original ordering. Our
objective is to find a permutation of $\mathcal{I}$ so that $K$  can be
approximated by an \hmatrix{}. The key is to define a distance  between a pair
of indices $i,j\in\mathcal{I}$, denoted as $d_{ij}$. Using the distances, we
then perform a hierarchical clustering of $\mathcal{I}$, which is used to
define the permutation, and determine which interaction go into
the sparse correction $S$ (using nearest neighbors). 

\paragraph{Three measures of distance.} We will define the point-based Euclidean distance, a Gram-space Euclidean distance, and a Gram-space angle distance.

\textbf{Geometric-$\ell^2$.} If we are given points $\{x_i\}_{i=1}^N$, then 
$d_{ij} = \lVert x_i - x_j \rVert_2$ is
the \emph{geometric $\ell^2$ distance}. This will be the \emph{geometry-aware}
reference implementation for cases where points are given.

\textbf{Gram-$\ell^2$ (or \emph{``kernel''} distance)}. Since $K$ is SPD, it is
the \emph{Gram matrix} of some set of \ipoint{ unknown Gram vectors},
$\{\phi_i\}_{i=1}^N \subset \mathbb{R}^N$  ( \cite{scholkopf-smola-02},
proposition 2.16, page 44).   That is, $K_{ij} = (\phi_i, \phi_j)$, where
$(\cdot,\cdot)$ denotes the $\ell^2$ inner product in $\mathbb{R}^N$.  Then we
define the \emph{Gram $\ell^2$ distance} as $d_{ij} = \lVert \phi_i - \phi_j
\rVert$. Computing the kernel distance only requires three entries of $K$:
\begin{equation}
  d_{ij}^2 = \lVert \phi_i \rVert^2 + \lVert \phi_j \rVert^2 - 2 (\phi_i, \phi_j) = K_{ii} + K_{jj} - 2K_{ij}. 
\end{equation}

\textbf{Gram angles (or \emph{``angle''} distance)}.
Our third measure of distance considers angles between Gram vectors, which is based on the standard sine distance (cosine similarity) in inner product spaces. 
We define the Gram angle distance as
$d_{ij} = \sin^2\left(\angle(\phi_i,\phi_j)\right)\in [0, 1]$.
This expression is chosen so that $d_{ij}$ is small for nearly collinear
Gram vectors, large for nearly orthogonal Gram vectors, and $d_{ij}$ is
inexpensive to compute. Although the value $d_{ij}$ may seem arbitrary,
we only compare values for the purpose of ordering,
so any equivalent metric will do.
Computing an angle distance only requires three entries of $K$:
\begin{equation}
  d_{ij} = 1 - \cos^2 \left(\angle(\phi_i,\phi_j)\right)
          = 1 - K_{ij}^2/(K_{ii} K_{jj}).
\end{equation}
To reiterate for emphasis, $d_{ij}$ define proper distances (metrics) because $K$ is SPD. And with distances, we can apply FMM.

\paragraph{Tree partitioning and nearest neighbor searches}
$K$ is permuted using a balanced binary tree.  The root node is assigned the full set of points, and the tree is constructed recursively by splitting a node's points evenly between two child nodes. The splitting terminates at nodes with some pre-determined \ipoint{leaf size} $m$. The leaf nodes then define a partial ordering of the indices: if leaf $\alpha$ is anywhere to the left of leaf $\beta$, then the indices of $\alpha$ precede those of $\beta$. We use this ordering to permute rows and columns of $K$. In the remainder of this paper, we use the notation $\alpha, \beta$ to refer interchangeably to a node or the set of indices belonging to the node.

 \begin{algorithm}[!t]
   \caption{{} $[\lc,\rc]=\texttt{metricSplit}(\alpha)$}
 \begin{algorithmic}
   \STATE $p = \mbox{argmax}\,( \{d_{ic}  \lvert i \in \alpha \})$; $q =
   \mbox{argmax}( \{ d_{ip} \lvert i \in \alpha \})$;
   \STATE $[\lc,\rc]=\mbox{medianSplit}(\{ d_{ip} - d_{iq} \lvert i \in
   \alpha\})$;
 \end{algorithmic}
 \label{a:split}
 \end{algorithm}  

In our implementation, we use a ball tree~\cite{march-xiao-yu-biros-sisc16}.
For geometric distances it costs $\bigO(N \log N)$. But Gram distances require
sampling to avoid $\bigO(N^2)$ costs.   Suppose we use one of the Gram
distances to split an interior node $\alpha$ between its left child $\lc$ and
right child $\rc$.  We define $c = \frac{1}{n_{\rm{c}}} \sum \phi_i$ to be an 
approximate centroid taken over a small sample of $n_{\rm{c}}$ Gram vectors belonging to $\alpha$. 
$n_{\rm{c}}$ is $\bigO(1)$.
Next, we find the point $p$ that is farthest away in distance from $c$, and the point $q$ that is farthest away from $p$. Then we split the indices $i \in \alpha$ on the values $d_{ip} - d_{iq}$, which measures the degree to which $i$ is closer to $p$ than to $q$. This approach is outlined in \algref{a:split}. 


We perform \emph{all nearest neighbors} (\texttt{ANN}) search  using randomized trees that are constructed in exactly the same way as the partitioning tree, except that $p$ and $q$ are chosen randomly. The search algorithm is described in~\cite{xiao-biros16} and (briefly) in the next section.


%% file: algfmm.tex
\begin{figure}[h]
  \centering
  \includegraphics[scale=.3]{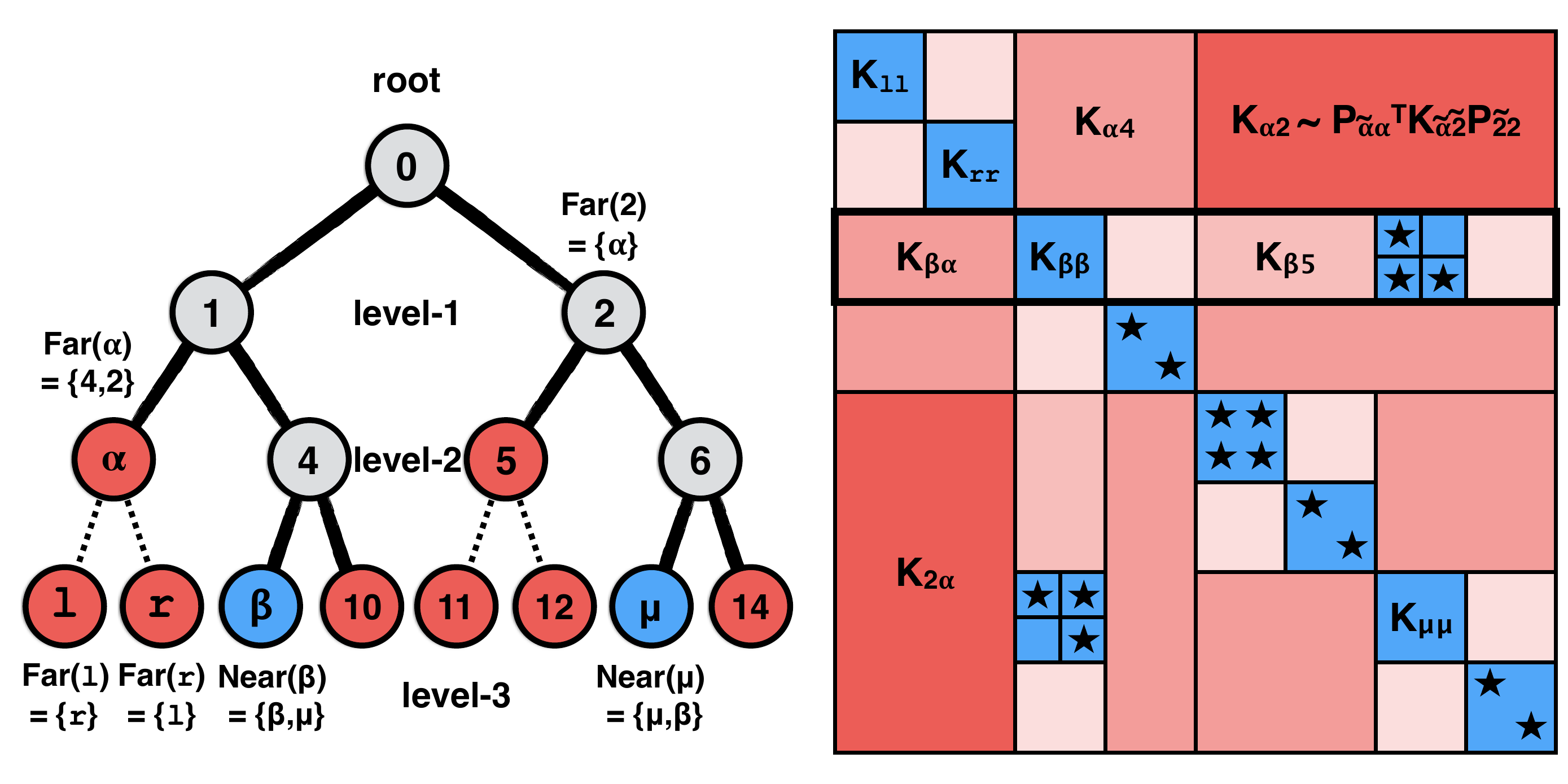}
  \caption{A partitioning tree (left) and corresponding hierarchically
    low-rank plus sparse matrix (right).
    The off-diagonal blocks are combinations of low-rank matrices (pink)
    and sparse matrices (blue).
    The $\bigstar$ symbol denotes an entry that cannot
    be approximated (because the corresponding interaction is between neighbors).
    The solid edges in the tree mark the path traversed by
    \texttt{FindFar}($\beta$,\texttt{0}).
    Since $K_{\beta\alpha}$ does not contain any neighbor interactions ($\bigstar$),
    this traversal adds $\alpha$ to \texttt{Far}$(\beta)$.
    In this example,
    \texttt{FindFar(\lc,0)} computes \texttt{Far}$(\lc)=\{\rc,4,2\}$, and
    \texttt{FindFar(\rc,0)} computes \texttt{Far}$(\rc)=\{\lc,4,2\}$.
    \algref{a:farfar} (\texttt{MergeFar}) then moves \texttt{Far}$(\lc) \cap Far(\rc)$
    into \texttt{Far}$(\alpha)$ so that
    \texttt{Far}$(\alpha)=\{4,2\}$, 
    \texttt{Far}$(\lc)=\{\rc\}$ and
    \texttt{Far}$(\rc)=\{\lc\}$.
  }
  \label{fig:tree}
\end{figure}

\hmatrix{} methods (including algebraic FMM) have two phases: \ipoint{compression} and \ipoint{evaluation}. 
As we discussed in the introduction, $K$ is compressed recursively using a binary tree
such that
\begin{equation}
  \label{e:partitioning}
  \sk{K}_{\alpha\alpha} =
\begin{bmatrix}
  \sk{K}_{\lc\lc} & 0 \\ 
  0 & \sk{K}_{\rc\rc} \\ 
\end{bmatrix} + 
\begin{bmatrix} 
0 & S_{\lc\rc} \\ 
S_{\rc\lc} & 0 \\ 
\end{bmatrix}+
\begin{bmatrix} 
  0 & UV_{\lc\rc} \\ 
  UV_{\rc\lc} & 0 \\ 
\end{bmatrix}, 
\end{equation} 
where $\lc$ and $\rc$ are \ipoint{left} and \ipoint{right child} of the treenode $\alpha$.
Each node $\alpha$ contains a set of matrix indices and the two 
children evenly split the indices such that $\alpha = \lc \cup \rc$.
(We overload the notation $\alpha$, $\beta$, $\lc$ and $\rc$ to
denote the matrix indices that those treenode own.)
In \figref{fig:tree}, the blue blocks depict $S$ (at all levels) and $D$ 
(in the leaf level), and the pink blocks depict the $UV$  matrices.

We use \ipoint{four tree traversals} to describe the algorithms in \gofmm{}:
postorder (\textbf{POST}), preorder (\textbf{PRE}),
any order (\textbf{ANY}),
and any order-leaves only (\textbf{LEAF}).
By  \ipoint{``task''} we refer to a computation that occurs when we visit a tree node during a traversal. 
We list all tasks required by \algref{a:compress} and \algref{a:evaluate} in \tabref{tab:tasks}.

We start by creating the binary metric ball tree in \algref{a:compress}
that represents the binary partitioning (and encodes a symmetric permutation of matrix $K$).
This requires the \emph{distance} metric $d_{ij}$ and a preorder traversal
(\textbf{PRE}) of the first task \texttt{SPLI($\alpha$)} in \tabref{tab:tasks}.

\begin{table}[!t]
\centering
{
\begin{tabular}{|>{\columncolor[gray]{0.8}}r|l|r|} 
\hline 
  \rowcolor[gray]{0.8}
Task & Operations & \texttt{FLOPS} \\
\hline 
\texttt{SPLI($\alpha$)} & split $\alpha$ into $\lc$ and $\rc$ \algref{a:split} & $\lvert \alpha \rvert$ \\
\hline 
\texttt{ANN($\alpha$)} & update $\MA{N}_{\alpha}$ with \texttt{KNN($K_{\alpha\alpha}$)} & $m^2$ \\ 
\hline 
\texttt{SKEL($\alpha$)} & $\sk{\alpha}$ in \algref{a:skeletonize} & $2s^3+2m^3$ \\ 
\hline 
\texttt{COEF($\alpha$)} & $P_{\sk{\alpha} \alpha}$ or $P_{\sk{\alpha} [\sk{\lc} \sk{\rc}]}$ in \algref{a:skeletonize} & $s^3$ \\ 
\hline 
\texttt{N2S($\alpha$)} & \texttt{\bf if} $\alpha$ is leaf \texttt{\bf then} $\sk{w}_{\alpha} = P_{\sk{\alpha}\alpha}w_{\alpha}$ & $2msr$ \\
                       & \texttt{\bf else} $\sk{w}_{\alpha} =
P_{\sk{\alpha}[\sk{\lc}\sk{\rc}]}[ \sk{w}_{\lc}; \sk{w}_{\rc}]$ & $2s^2r$\\
\hline
\texttt{SKba($\beta$)}  & $\forall\alpha\in\mathtt{Far}(\beta)$,
$K_{\sk{\beta}\sk{\alpha}} = K(\sk{\beta}, \sk{\alpha})$ & $ds^2\lvert \mathtt{Far}(\beta) \rvert$ \\ 
\hline 
\texttt{S2S($\beta$)} & $\sk{u}_{\beta} =
\sum_{\alpha\in \mathtt{Far}(\beta)} K_{\sk{\beta}\sk{\alpha}}\sk{w}_{\alpha}$
& $2s^2r\lvert \mathtt{Far}(\beta) \rvert$\\ 
\hline 
\texttt{S2N($\beta$)} & \texttt{\bf if} $\alpha$ is leaf \texttt{\bf then} $u_{\beta} = P_{\sk{\beta}\beta}^{T}\sk{u}_{\beta}$ & $2msr$ \\ 
                      & \texttt{\bf else} $[\sk{u}_{\lc};\sk{u}_{\rc}] += P_{\sk{\beta}[\sk{\lc}\sk{\rc}]}^{T}\sk{u}_{\beta}$ & $2s^2r$ \\
\hline
\texttt{Kba($\beta$)}  & $\forall\alpha\in\mathtt{Near}(\beta)$, $K_{\beta\alpha} = K(\beta, \alpha)$ & $m^2\lvert\mathtt{Near}(\beta)\rvert$ \\ 
\hline  
\texttt{L2L($\beta$)} & $u_{\beta} += \sum_{\alpha\in
\mathtt{Near}(\beta)}K_{\beta\alpha}w_{\alpha}$ & $2m^2r\lvert\mathtt{Near}(\beta) \rvert$\\ 
\hline 
\end{tabular}
}
\caption{Tasks and their costs in \texttt{FLOPS}.
  \texttt{SPLI} (tree splitting), \texttt{ANN} (all nearest-neighbors), \texttt{SKEL} (skeletonization), 
  \texttt{COEF} (interpolation) \texttt{SKba} and 
  \texttt{Kba} (caching submatrices) occur in the compression phase.
  Interactions \texttt{N2S} (nodes to skeletons), \texttt{S2S} (skeletons to skeletons),
  \texttt{S2N} (skeletons to nodes), and \texttt{L2L} (leaves to leaves) occur in 
  the evaluation phase.}
\label{tab:tasks}
\end{table} 

\begin{algorithm}[!t]
\caption{{} \texttt{Compress}($K$)}
\begin{algorithmic}[1]
  \STATE \texttt{\bf for each} randomized tree \texttt{\bf do} \# iterative neighbor search
  \STATE \gap (\texttt{\bf PRE}) \texttt{SPLI($\alpha$)} \# create a random projection tree 
  \STATE \gap (\texttt{\bf LEAF}) \texttt{ANN($\alpha$)} \# search $\kappa$ neighbors in leaf nodes
  \STATE (\texttt{\bf PRE}) \texttt{SPLI($\alpha$)} \# create a metric ball tree 
  \STATE (\texttt{\bf LEAF}) \texttt{LeafNear($\beta$)} \# build \texttt{Near}$(\beta)$ using $\MA{N}(\beta)$
  \STATE (\texttt{\bf LEAF}) \texttt{FindFar($\beta$,root)} \# find \texttt{Far}$(\beta)$ using \texttt{MortonID}
  \STATE (\texttt{\bf POST}) \texttt{MergeFar($\alpha$)} \# merge \texttt{Far}$(\lc)$, \texttt{Far}$(\rc)$ to \texttt{Far}$(\alpha)$
  \STATE (\texttt{\bf POST}) \texttt{SKEL($\alpha$)} \# compute skeletons $\sk{\alpha}$
  \STATE (\texttt{\bf ANY}) \texttt{COEF($\alpha$)} \# compute the coefficient matrix $P$
  \STATE (\texttt{\bf ANY}) \texttt{Kba($\beta$)} \# optionally evaluate and cache $K_{\beta\alpha}$ 
  \STATE (\texttt{\bf ANY}) \texttt{SKba($\beta$)} \# optionally evaluate and cache $K_{\sk{\beta}\sk{\alpha}}$
\end{algorithmic}
\label{a:compress}
\end{algorithm}

\begin{algorithm}[!t]
\caption{{} \texttt{LeafNear}($\beta$)}
\begin{algorithmic}
  \STATE  \texttt{Near}$(\beta)$ = \{\texttt{MortonID}$(i)$ : $\forall i\ \in \MA{N}(\beta)$\}
\end{algorithmic}
\label{a:nearnear}
\end{algorithm}

\begin{algorithm}[!t]
  \caption{{} \texttt{FindFar}($\beta=\mathtt{leaf}$, $\alpha$)}
\begin{algorithmic}
  \STATE \texttt{\bf if} $\alpha \cap Near(\beta) \neq \phi$ using \texttt{MortonID} \texttt{\bf then}
  \STATE \gap \texttt{FindFar}($\beta$,\lc); \texttt{FindFar}($\beta$,\rc); 
  \STATE \texttt{\bf else} \texttt{Far}$(\beta) = Far(\beta) \cup \alpha$;
\end{algorithmic}
\label{a:nearfar}
\end{algorithm}

\begin{algorithm}[!t]
\caption{{} \texttt{MergeFar}($\alpha$)}
\begin{algorithmic}
  \STATE \texttt{MergeFar}($\lc$); \texttt{MergeFar}($\rc$);
  \STATE \texttt{Far}$(\alpha) = Far(\lc) \cap Far(\rc)$; 
  \STATE \texttt{Far}$(\lc) = Far(\lc) \backslash Far(\alpha)$; \texttt{Far}$(\rc) = Far(\rc) \backslash Far(\alpha)$;
\end{algorithmic}
\label{a:farfar}
\end{algorithm}

\textbf{Node lists and near-far pruning.}
\gofmm{} tasks require that every tree node maintains three lists.
For a node $\alpha$, these lists are the neighbor list $\MA{N}(\alpha)$,
near interaction list $\mathtt{Near}(\alpha)$, and far interation list
$\mathtt{Far}(\alpha)$. Computing these lists requires
defining neighbors for indices based on the distance $d_{ij}$ 
and the Morton ID.


A pair of nodes $\alpha$ and $\beta$ is said to be \emph{far} if
$K_{\beta\alpha}$ is low-rank and \emph{near} otherwise.
We use neighbor-based pruning~\cite{march-xiao-yu-biros-sisc16}
to determine the near-far relation.
Neighbors are  defined based on the specified distance
$d_{ij}$. For each $i$, we search for the $\kappa$ indices $j$
that result in the smallest $d_{ij}$.
The Morton ID is a bit array that codes the path from the root to a tree node or index $i$. The Morton ID of an index $i$ is the Morton ID of the leaf node (in \gofmm{} ball metric tree) that contains it.
We use \texttt{MortonID()} to denote this.

\ipoint{Index nearest neighbor list} $\MA{N}$($i$): As we discussed, \gofmm{}
requires a preprocessing step in which we compute the nearest neighbors for 
every index $i$ using a greedy search (steps 1--3 in~\algref{a:compress}). 
This constructs a list of $\kappa$ nearest-neighbor for each $i\in\alpha$
iteratively.
In each iteration, we create a randomized projection
tree~\cite{dasgupta-freund08,liberty2007randomized,march-xiao-yu-biros-sisc16},
and we search for neighbors of $i$ only in the leaf node $\alpha$ that contains
$i$ using an exhaustive search~\cite{yu2015performance}. That is, for each $i\in\alpha$, 
we only search for small $d_{ij}$ where $j\in\alpha$ as well. 
The iteration stops after reaching $80\%$ accuracy or 10 iterations.

\ipoint{Node neighbor list} $\MA{N}$($\alpha$): 
Then we construct the neighbor list $\MA{N}(\alpha)$ of a leaf node 
$\alpha$ by merging all neighbors of $i \in \alpha$. For non-leaf nodes 
the list is constructed recursively~\cite{march-xiao-biros-fmm-e15}.

\ipoint{Near list of a node} \texttt{Near}$(\alpha)$:
Leaf nodes $\alpha, \beta$ are considered near if 
$\alpha \cap \MA{N}(\beta)$ is nonempty
(i.e., $K_{\alpha \beta}$ contains at least one neighbor 
($\bigstar$) in \figref{fig:tree}). The \texttt{Near} list is defined only for leaf nodes and contains only leaf nodes.  For each leaf node $\beta$, \texttt{Near}$(\beta)$ is constructed
using \texttt{LeafNear} (\algref{a:nearnear}). For each neighbor $i\in\MA{N}(\beta)$, \texttt{LeafNear($\beta$)}  adds \texttt{MortonID}$(i)$ to \texttt{Near}$(\beta)$. 
Notice that the size of \texttt{Near}$(\beta)$ determines the number
of direct evaluations in the off-diagonal blocks. To prevent the cost from growing too fast,
we introduce a user-defined parameter \textbf{budget} such that
\begin{equation}
  \lvert \mathtt{Near}(\beta)\rvert < \mathrm{budget} \times (N/m). 
\label{e:budget}
\end{equation}
While looping over $i\in\MA{N}(\beta)$, instead of directly adding \texttt{MortonID}$(i)$
to \texttt{Near}$(\beta)$, we only mark it with a ballot.
Then we insert candidates to \texttt{Near}$(\beta)$ according to their votes
until \eqref{e:budget} is reached.
To enforce symmetry of $\tilde{K}$, we loop over all \texttt{Near} lists and enforce the following:
if $\alpha \in$ \texttt{Near}$(\beta)$ then $\beta \in$ \texttt{Near}$(\alpha)$.


\ipoint{Far list of a node} \texttt{Far}$(\alpha)$:
\texttt{Far}$(\alpha)$ is constructed in two steps in \algref{a:compress}.
First for each leaf node $\beta$, we invoke \texttt{FindFar}($\beta$, \texttt{root}) (\algref{a:nearfar}).
Upon visiting $\alpha$, we check whether $\alpha$ is a parent of any 
leaf node in \texttt{Near}$(\beta)$ using \texttt{MortonID}.
If so, we add $\alpha$ to \texttt{Far}$(\beta)$; otherwise, we recurse to the 
two children of $\alpha$.
The second step is a postorder traversal on \texttt{MergeFar(root)}
(\algref{a:farfar}).
This process merges the common nodes from two children lists \texttt{Far}$(\lc)$ and \texttt{Far}$(\rc)$.
These common nodes are removed from the children and added to 
their parent list \texttt{Far}$(\alpha)$.
In \figref{fig:tree}, 
\texttt{FindFar} can be identified by the smallest square pink blocks, 
and \texttt{MergeFar} merges small pink blocks into larger blocks.

\textbf{Low-rank approximation.}
We approximate off-diagonal matrix blocks with a nested
interpolative decomposition (ID)~\cite{halko-martinsson-tropp11}.
Let $\beta$ be the indices in a leaf node and
$I=\{1,...,N\}\backslash \beta$ be the set complement.
The \emph{skeletonization} of $\beta$ is a rank-$s$
approximation of its off-diagonal blocks $K_{I\beta}$ using the ID, 
which we write as
\begin{equation}
K_{I\beta} \approx K_{I\sk{\beta}}P_{\sk{\beta}\beta},
\end{equation}
where $\sk{\beta} \subset \beta$ is the \emph{skeleton} of $\beta$.
$K_{I\sk{\beta}} \in \mathbb{R}^{(N - \lvert \beta \rvert) \times s}$ is a column
submatrix of $K_{I \beta}$, and $P_{\sk{\beta}\beta} \in
\mathbb{R}^{s\times \lvert \beta \rvert}$ is a matrix of interpolation
coefficients, where $\ns$ is the approximation rank.

To efficiently compute this approximation, we select a sample subset $I'
\subset I$
using neighbor-based importance sampling~\cite{march-xiao-yu-biros-sisc16}. 
We then perform a rank-revealing QR factorization (\texttt{GEQP3}) on $K_{I' \beta}$.
The skeletons $\sk{\beta}$ are selected to be the first $\ns$ pivots, and
the matrix $P_{\sk{\beta}\beta}$ is computed by a triangular solve (\texttt{TRSM})
using the triangular factor $R$.
The rank $\ns$ is chosen adaptively such that $\sigma_{\ns+1}(K_{I' \beta}) < \idtol$,
where $\sigma_{\ns+1}(K_{I' \beta})$ is the estimated $\ns+1$ singular value and
$\idtol$ is related to a user-specified error tolerance.

For an internal node $\alpha$, we form the skeletonization in the same way,
except that the columns are also sampled using the skeletons of the children of $\alpha$.
That is, the ID is computed for $K_{I'[\sk{\lc}\sk{\rc}]}$,
where $[\sk{\lc}\sk{\rc}] = \sk{\lc} \cup \sk{\rc}$ contains the skeletons of the children
of $\alpha$:
\begin{equation}
K_{I[\sk{\lc}\sk{\rc}]} \approx K_{I\sk{\alpha}}P_{\sk{\alpha}[\sk{\lc}\sk{\rc}]}.
\end{equation}
This way, the skeletons are nested:
$\sk{\alpha} \subset \sk{\lc} \cup \sk{\rc}$.

As a consequence of the nesting property,
we can use $P_{\sk{\lc}\lc}$ and $P_{\sk{\rc}\rc}$ to construct an approximation of
the full block $K_{I\alpha}$:
\begin{equation}
K_{I\alpha}
\approx
K_{I[\sk{\lc}\sk{\rc}]}
\begin{bmatrix}
  P_{\sk{\lc}\lc} & \\
                  & P_{\sk{\rc}\rc} \\
\end{bmatrix}
\approx
K_{I\sk{\alpha}} P_{\sk{\alpha}[\sk{\lc}\sk{\rc}]}
\begin{bmatrix}
  P_{\sk{\lc}\lc} & \\
                  & P_{\sk{\rc}\rc} \\
\end{bmatrix}.
\end{equation}
Then we have a \emph{telescoping} expression for the full coefficient matrix:
\begin{equation}
P_{\sk{\alpha}\alpha} =
P_{\sk{\alpha}[\sk{\lc}\sk{\rc}]}
\begin{bmatrix}
  P_{\sk{\lc}\lc} & \\
                  & P_{\sk{\rc}\rc} \\
\end{bmatrix}.
\label{e:telescope}
\end{equation}
We never explicitly form $P_{\sk{\alpha}\alpha}$, but instead use the
telescoping expression during evaluation.

\begin{algorithm}[!t]
\caption{{} [$\sk{\alpha}, P_{\sk{\alpha} \alpha}$]=\texttt{Skeleton}($\alpha$)}
\begin{algorithmic}
  \STATE \texttt{\bf if} $\alpha$ is leaf \texttt{\bf then return} [$\sk{\alpha}, P_{\sk{\alpha} \alpha}$] = {\tt ID}$(\alpha)$;
  \STATE $[\sk{\lc},]=\texttt{Skeleton}(\lc)$; $[\sk{\rc},]=\texttt{Skeleton}(\rc)$;
  \STATE \texttt{\bf return} [$\sk{\alpha}, P_{\sk{\alpha} [\sk{\lc} \sk{\rc}]
  }$] = {\tt ID}($[\sk{\lc} \sk{\rc}]$);
\end{algorithmic}
\label{a:skeletonize}
\end{algorithm}

\algref{a:skeletonize} computes the skeletonization for all tree nodes with a 
postorder traversal.
There are two tasks for each tree node $\alpha$ listed in \tabref{tab:tasks}: 
(1) \texttt{SKEL($\alpha$)} selects $\sk{\alpha}$ (in the critical path) and
(2) \texttt{COEF($\alpha$)} computes $P_{\sk{\alpha}[\sk{\lc}\sk{\rc}]}$.
Notice that in \algref{a:compress} only \texttt{SKEL($\alpha$)} needs to be
executed in postorder (\textbf{POST}), but \texttt{COEF($\alpha$)} can be in any
order (\textbf{ANY})
as long as \texttt{SKEL($\alpha$)} is finished. Such parallelism can
only be specified at the task level, which later inspires our task-based
parallelism in \secref{s:par}.
At the end of the compression, we can optionally evaluate and cache 
all $K_{\beta\alpha}$ in \texttt{Near}$(\beta)$ and all
$K_{\sk{\beta}\sk{\alpha}}$ in \texttt{Far}$(\beta)$ by executing \texttt{Kba($\beta$)} 
and \texttt{SKba($\beta$)} in any order.
Given enough memory (at least $\MA{O}(N)$ for all 
$K_{\beta\alpha}$ and $\K_{\sk{\beta}\sk{\alpha}}$), caching 
can reduce the time spent evaluating and gathering submatrices. 

\textbf{Evaluation.}
Following~\cite{march-xiao-biros-fmm-e15}, we present
\algref{a:evaluate} a four-step process for computing
\eqref{e:appromatvec}.
The idea is to approximate each \textbf{matvec} $u_{\beta}+=K_{\beta\alpha}w_{\alpha}$
in \texttt{Far}$(\beta)$ using 
a two-sided ID to accumulate
$P_{\sk{\beta}\beta}^{T}K_{\sk{\beta}\sk{\alpha}}P_{\sk{\alpha}\alpha}w_{\alpha}$,
where $P_{\sk{\alpha}\alpha}, P_{\sk{\beta}\beta}$ are given by the
telescoping expression \eqref{e:telescope}. 
For more details, see~\cite{march-xiao-biros-fmm-e15}.

\begin{algorithm}
\caption{{} \texttt{Evaluate}($u$,$w$)}
\begin{algorithmic}[1]
  \STATE (\textbf{POST}) \texttt{N2S($\alpha$)} \# compute skeleton weights $\sk{w}$
  \STATE (\textbf{ANY}) \texttt{S2S($\beta$)} \# apply skeleton basis $K_{\sk{\beta}\sk{\alpha}}$
  \STATE (\textbf{PRE}) \texttt{S2N($\beta$)} \# accumulate skeleton potentials $\sk{u}$
  \STATE (\textbf{ANY}) \texttt{L2L($\beta$)} \# accumulate direct \textbf{matvec} to $u$
\end{algorithmic}
\label{a:evaluate}
\end{algorithm}

The first step is to perform a postorder traversal (\textbf{POST}) on 
\texttt{N2S($\alpha$)} (Nodes To Skeletons).
This computes the \emph{skeleton weights} $\sk{w}_{\alpha}=P_{\sk{\alpha}\alpha}w_{\alpha}$
for each leaf node, and 
$\sk{w}_{\alpha} = P_{\sk{\alpha}[\sk{\lc}\sk{\rc}]}[ \sk{w}_{\lc}; \sk{w}_{\rc}]$
for each inner node.
Recall that in \texttt{COEF($\alpha$)}, we have computed $P_{\sk{\alpha}\alpha}$
for each leaf node and $P_{\sk{\alpha}[\sk{\lc}\sk{\rc}]}$ for each 
internal node. 
\texttt{S2S($\beta$)} (Skeletons to Skeletons) applies the \emph{skeleton basis} $K_{\sk{\beta}\sk{\alpha}}$
and accumulates \emph{skeleton potentials} $\sk{u}$ for each node:
$\sk{u}_{\beta}=\sum_{\alpha\in Far(\beta)}
K_{\sk{\beta}\sk{\alpha}}\sk{w}_{\alpha}$.
As soon as $\sk{w}_{\alpha}$ are computed in \texttt{N2S}, \texttt{S2S} can
be executed in any order.
\texttt{S2N($\beta$)} (Skeletons To Nodes) performs interpolation on the left
and accumulates
$\sk{u}$ with a preorder traversal. This uses the transpose of
\eqref{e:telescope}.
For each node $\beta$, we accumulate
$[\sk{u}_{\lc};\sk{u}_{\rc}] +=
P_{\sk{\beta}[\sk{\lc}\sk{\rc}]}^{T}\sk{u}_{\beta}$ to its children.
In the leaf node, $u_{\beta} = P_{\sk{\beta}\beta}^{T}\sk{u}_{\beta}$
directly accumulates to the output.
These three tasks compute all \textbf{matvec} for the \emph{far} nodes (pink
blocks in \figref{fig:tree}).
All \textbf{matvec} on $K_{\beta\alpha}$ in \texttt{Near}$(\beta)$ (blue blocks)
are computed by \texttt{L2L($\beta$)} (Leaves To Leaves) and directly accumulated to $u_{\beta}$.

\textbf{Complexity.}
The worst case cost of \algref{a:evaluate} is $\MA{O}(N^2)$, when 
$\lvert\mathtt{Near}(\alpha)\rvert=(N/m)$ for all $\alpha$.
The best case occurs when each \texttt{Near}$(\alpha)$ only contains $\alpha$ itself.
We fix the rank $s$ and 
leaf size $m$.
The tree has $\MA{O}(N/m)$ leaf nodes and $\MA{O}(N/m)$ interior nodes, so in the best case,
overall
\texttt{N2S} has $\MA{O}(2ms(N/m)+2s^2(N/m))$ work,
\texttt{S2S} has  $\MA{O}(2s^2(N/m))$ work,
\texttt{S2N} has $\MA{O(2ms(N/m)+2s^2(N/m))}$ work,
and \texttt{L2} has $\MA{O}(2m^2(N/m))$.
When $s$ and $m$ are held constant, the total work is $O(N)$ per right hand side.
In \gofmm{}, this is controlled by the \textbf{budget}.

%% file: par.tex
In \hmatrix{} methods and FMM the main algorithmic pattern is a tree traversal. A traversal may exhibit high parallelism at the leaf level but due to the dependencies the parallelism typically diminishes near  the root level. In addition, if the workload per tree node varies, load balancing becomes an issue.  Most static scheduling codes employ level-by-level traversals, which introduces unnecessary synchronizations. In \gofmm{}, we observe significant workload variations during the compression (\algref{a:skeletonize}) and during the evaluation (tasks \texttt{N2S} and \texttt{S2N}). 

One solution is to exploit parallelism in finer granularity. For example, when the number of tree nodes
in the single tree level is less than the number of cores, we can use
multi-threaded BLAS/LAPACK on a single tree node. However, this is insufficient
if the workload does not increase significantly (e.g. growing with
$\lvert\alpha\rvert$)
while approaching the root. 
(That is, the workload must be within the strong scaling range of BLAS/LAPACK to
be efficient).

To partially address these challenges we abandon the convenient level-by-level
traversal and explore an \emph{out-of-order} approach using dynamic scheduling.
To this end we test two approaches and compare them with a level-by-level
traversal. In the first approach, we introduce a self-contained runtime system.
In the second approach we test the same ideas with \texttt{OpenMP}'s \texttt{omp task depend} feature. 

\begin{figure}[!t]
  \centering
  \includegraphics[scale=.4]{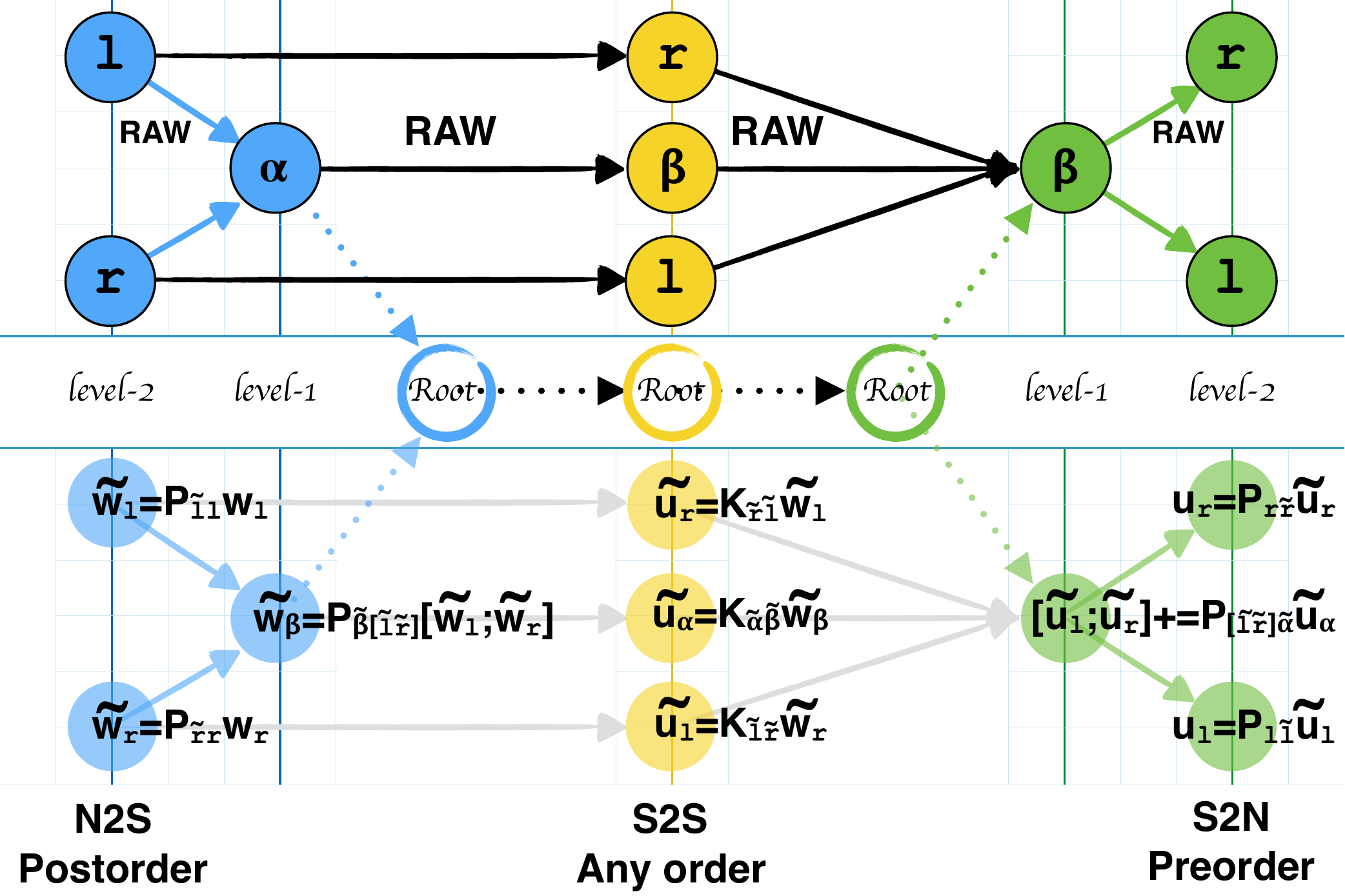}
  \caption{Dependency graph for steps 1--3 of \algref{a:evaluate} (step 4 is completely independent of steps 1--3).
    Each tree node denotes a task, and the arrows between nodes imply a
    dependency. Here $\mathtt{Near}(\alpha)$ only contains itself (HSS). For example, yellow node $\beta$ has a \textbf{RAW}
    dependency following blue $\alpha$, because \texttt{S2S($\beta$)} computes 
  $\sk{u}_{\beta} = \sum_{\alpha \in
  \mathtt{Near}(\beta)}K_{\sk{\beta}\sk{\alpha}}\sk{w}_{\alpha}$.
  When $\mathtt{Near}(\beta)$ contains more than just itself. The dependencies
    are unknown at compile time and thus, \texttt{omp task depend} fails to describe 
    the dependencies between \texttt{N2S} and \texttt{S2S}.
  }
  \label{fig:n2s2n}
\end{figure}

\textbf{Dependency analysis.}
Recursive preorder and postorder traversals 
inherently encode \textbf{R}ead/\textbf{W}rite dependencies 
between tree nodes. Following \algref{a:compress} and \algref{a:evaluate}, we can describe dependencies between different tasks. However, due to dynamic granularity of tasks we need a data flow 
analysis at runtime.  For example, dependencies between \texttt{N2S} and \texttt{S2S} 
cannot be discovered at compile time, because the \textbf{RAW}
(read after write) dependencies on $\sk{w}_{\alpha}$ are computed by neighbors
$\MA{N}(\alpha)$. In order to build dependencies at runtime as a direct acyclic graph (DAG), 
we perform a \emph{symbolic} execution on \algref{a:compress} 
and \algref{a:evaluate}. For simplicity, below we just discuss the evaluation phase for the HSS case 
(the FMM case is more involved).

\figref{fig:n2s2n} depicts task dependencies (by tasks we mean algorithmic
tasks defined in~\tabref{tab:tasks}) during the evaluation
phase~\algref{a:evaluate} for \texttt{N2S}, \texttt{S2S} and \texttt{S2N} where
the off-diagonal blocks are  low-rank (HSS) with $S=0$. This task dependency
tree is generated by our runtime using symbolic traversals. The \texttt{N2S},
\texttt{S2S}, and \texttt{S2N} execution order is performed on a binary
tree\footnote{\scriptsize Execution order from left to right: dependencies are
easier to follow if one rotates the page by $90^{\circ}$ counter-clockwise}.

We use three symbolic traversals of~\algref{a:evaluate}. In the first traversal (postorder) we
find that  $\sk{w}_{\lc}$ is written by $\lc$. Going from $\sk{w}_{\lc}$ to
$\sk{w}_\beta$, we annotate that $\sk{w}_{\lc}$ is read by $\beta$, i.e.
$\sk{w}_{\beta}= P_{\sk{\beta}[\sk{\lc}\sk{\rc}]}[\sk{w}_{\lc};\sk{w}_{\rc}]$.
This \textbf{RAW} dependency is an edge from $\lc$ to $\beta$ in the DAG.

Intertask dependencies are discovered by  the \emph{symbolic}
execution of the yellow tree. At node $\beta$ (in yellow), the relation
$\sk{u}_{\alpha}=K_{\sk{\alpha}\sk{\beta}\sk{w}_{\beta}}$
will read $\sk{w_{\beta}}$. Again this is a \textbf{RAW} dependency, hence the edge from the blue $\beta$ to the yellow $\alpha$. The whole dependency graph for steps 1--3 is built after the green postorder tree traversal. Step 4 in \algref{a:evaluate} is independent of steps 1--3. Although this runtime data flow analysis has some overhead, the amount is almost negligible ($<1\%$) compared to the total execution time.  

\textbf{Runtime.}
With a dependency graph, scheduling can be done in 
\emph{static} or \emph{dynamic} fashion. 
Due to unknown adaptive rank $s$ at compile time, we implement a 
light-weight dynamic Heterogeneous Earliest Finish Time 
(HEFT)~\cite{topcuoglu2002performance} using \texttt{OpenMP} threads.
Each worker (thread) in the runtime can use more than one physical 
core with either  a nested \texttt{OpenMP} construct or by employing a device  (accelerator) as a slave.
Tasks that satisfy all dependencies in the DAG  will be 
dispatched to a ``ready'' queue. Each worker keeps 
consuming tasks in its own queue until no tasks are left.

Although we estimate a cost for each task\footnote{\scriptsize 
We divide costs for tasks by the theoretical peak 
\texttt{FLOPS} of the target architecture and a discount factor. 
For memory-bound tasks we use the theoretical \texttt{MOPS}
instead.} 
in \tabref{tab:tasks}, the runtime of a normal worker (or one with an accelerator) depends on the problem and can only be determined at runtime.
The HEFT schedule is implemented using an estimated finish time of all pending tasks in a specific worker's  ready queue. Each task dispatched from the DAG is assigned to a queue such that the maximum estimated finish time of each queue  is minimized. For the case where the estimation is inaccurate, we also implement
a job stealing mechanism. 

\textbf{Other parallel implementation.}
We briefly introduce other possible parallel implementations and conduct a strong scaling experiment in \secref{s:results}.
Here we implemented parallel level-by-level traversals for all tasks 
that require preorder and postorder traversals and do not exploit out-of-order parallelism. 
For tasks that can be executed in any order, we simply use 
\texttt{omp parallel for} with dynamic scheduling. 
If there are not enough tree nodes in a tree level, we use
nested parallelism with inner \texttt{OpenMP} constructs 
and multi-threaded BLAS/LAPACK.

The \texttt{omp task} version is implemented using recursive 
preorder or postorder traversals. Due to the overhead 
of the deep call stack, this implementation can be much slower than others. Although we tested it we do not report results because it is not competitive.

We also implemented (and report results for) \texttt{omp task depend}, since
\texttt{OpenMP-4.5} supports task parallelism with dependencies. However there
are two issues. First,
\texttt{omp task depend} requires all dependencies to be known at compile time,
which is not the case for the FMM (tasks \texttt{N2S} and \texttt{S2S}).
Second, without knowledge of the estimated finish time, the \texttt{OpenMP} scheduler will be suboptimal.  Finally for CPU-GPU hybrid architectures, scheduling GPU tasks purely with 
\texttt{omp task} can be very challenging. 

\textbf{CPU-GPU hybrid.}
GPUs usually offer high computing capacity, but  performance 
can easily be bounded by the PCI-E bandwidth. 
Because most computations in 
\algref{a:compress} are complex and memory bound\footnote{\scriptsize Although \texttt{GEQP3} 
and \texttt{TRSM} can be performed on GPUs with \texttt{MAGMA}
(\url{http://icl.cs.utk.edu/magma/})
and \texttt{cublas}, we find this inefficient for our methods.},
we do not use GPUs for the compression.
Instead we only pre-fetch submatrices $K_{\beta\alpha}$ and
$K_{\sk{\beta}\sk{\alpha}}$ to the device memory to overlap with
computations on the host (CPUs).
During the evaluation, our runtime will decide--depending on the number of
\texttt{FLOPS}-- whether to issue a batch
of tasks (up to 8) to the GPU in concurrent (using \texttt{stream}). This usually 
occurs in \texttt{N2S} and \texttt{S2N} where the size of 
\texttt{cublasXgemm}
is bounded by $s$ and $m$.
Furthermore, to hide communication time between CPU and GPU, 
all arguments of the next task in queue are pre-fetched using
asynchronous communication for pipelining. 
Finally, because a worker with a GPU is usually 
$50\times$ to $100\times$ more capable than others, we disable job stealing balancing for GPU workers.
This optimization prevents the GPU from idling.


%% file: setup.tex
We perform  experiments on Haswell, KNL, ARM, and NVIDIA GPU architectures
with four different setups to examine the accuracy and efficiency of
our methods. We demonstrate 
(1) the robustness and effectiveness of our geometry-oblivious FMM,
(2) the scalability of our runtime system against other parallel schemes, 
(3) the accuracy and cost comparison with other software, and 
(4) the absolute efficiency (in percentage of peak performance).

\textbf{Implementation and hardware.} Please refer to \secref{s:sup_setup} 
for all configuration in the reproducibility artifact section. Our
tests were conducted on TACC's Lonestar 5, (two 12-core, 2.6GHz, Xeon
E5-2690 v3 ``Haswell''), TACC's Stampede 2 (68-core, 1.4GHz, Xeon Phi
7250 ``KNL'') and CSCS's Piz Daint (12-core, 2.3GHz, Xeon E5-2650 v3
and NVIDIA Tesla P100).

\paragraph{\textbf{Matrices}}
We generated 22 matrices emulating different problems.
{\bf K02} is a 2D regularized inverse Laplacian squared, 
resembling the Hessian operator of a PDE-constrained optimization problem. 
The Laplacian is discretized using a 5-stencil finite-difference scheme 
with Dirichlet boundary conditions on a regular grid. 
{\bf K03} has the same setup with the 
oscillatory Helmholtz operator and 10 points per wave length. 
{\bf K04--K10} are kernel matrices in six dimensions 
(Gaussians with different bandwidths, narrow and wide; 
Laplacian Green's function, polynomial and cosine-similarity).
{\bf K12--K14} are 2D advection-diffusion operators on a regular 
grid with highly variable coefficients. 
{\bf K15,K16} are 2D pseudo-spectral  advection-diffusion-reaction 
operators with variable coefficients. 
{\bf K17} is a 3D pseudo-spectral operator with variable coefficients. 
{\bf K18} is the inverse squared Laplacian in 3D with variable coefficients. 
{\bf G01--G05} are the inverse Laplacian of the  \textbf{powersim}, 
\textbf{poli\_large}, \textbf{rgg\_n\_2\_16\_s0}, 
\textbf{denormal}, and \textbf{conf6\_0-8x8-30} graphs from 
\href{http://yifanhu.net/GALLERY/GRAPHS/search.html}{UFL}.  

\textbf{K02}--\textbf{K03}, 
\textbf{K12}--\textbf{K14}, and \textbf{K18} resemble inverse 
covariance matrices and Hessian operators from optimization and 
uncertainty quantification problems. 
\textbf{K04}--\textbf{K10} resemble 
classical kernel/Green function matrices but in high dimensions. 
\textbf{K15}--\textbf{K17} resemble pseudo-spectral operators.
\textbf{G01}--\textbf{G05} ($N=$ 15838, 15575, 65536, 89400, 49152)
are graphs for which we do not have geometric information.
For \textbf{K02}--\textbf{K18}, we use $N=65536$ if not specified.

Also, we use kernel matrices from machine learning: 
\textbf{COVTYPE} (100K, 54D, cartographic variables); and
\textbf{HIGGS} (500K, 28D,  physics)~\cite{Lichman:2013};
\textbf{MNIST} (60K, 780D, digit recognition)~\cite{chang2011libsvm}.
For these datasets, we use a Gaussian kernel with bandwidth $h$.

\gofmm{} supports both double and single precision. All experiments with matrices \textbf{K02--K18} and \textbf{G01--G05} are in single precision. The results for {\bf COVTYPE, HIGGS, MNIST} are in double precision.

\textbf{Parameter selection and accuracy metrics.}
We control $m$ (leaf node size), 
$s$ (maximum rank), $\tau$ (adaptive tolerance), $\kappa$ (number of
neighbors), \ipoint{budget} (a key paramter for amount of direct evaluations
and for switching between HSS and FMM) and partitioning (\textbf{Kernel},
\textbf{Angle}, \textbf{Lexicographic}, \textbf{geometric}, \textbf{random}).
We use $m=$256--512; on average this gives good overall time.  The adaptive
tolerance $\tau$, reflects the error of the subsampled  block and may not
correspond to the output error $\epsilon_2$. Depending on the problem, $\tau$ may misestimate the rank.
Similarly, this may occur in \texttt{HODLR}, \texttt{STRUMPACK} and \texttt{ASKIT}.
We use $\tau$ between 1E-2 and 1E-7, $s=m$, $k=32$ and $3\%$ budget. To enforce a HSS approximation, we use $0\%$ budget. The Gaussian bandwidth values are taken from~\cite{march-xiao-biros-e15} and produce optimal learning rates.



Throughout we use relative error $\epsilon_2$ defined as the following
\begin{equation}
  \epsilon_2 = \|\sk{K}w-Kw\|_F / \|Kw\|_F,\ \mbox{where} \ w \in \mathbb{R}^{N\times r}.
\end{equation}
This metric requires $\MA{O}(rN^2)$ work; to reduce the computational effort we instead sample 100 rows of $K$.
In all tables, we use ``\rm{Comp}'' and ``\rm{Eval}'' to refer the the 
compression and evaluation time in seconds, and ``\rm{GFs}'' to 
\rm{GFLOPS} per node.


%% file: results.tex
\begin{figure}[!t]
  \centering
  \includegraphics[scale=.3]{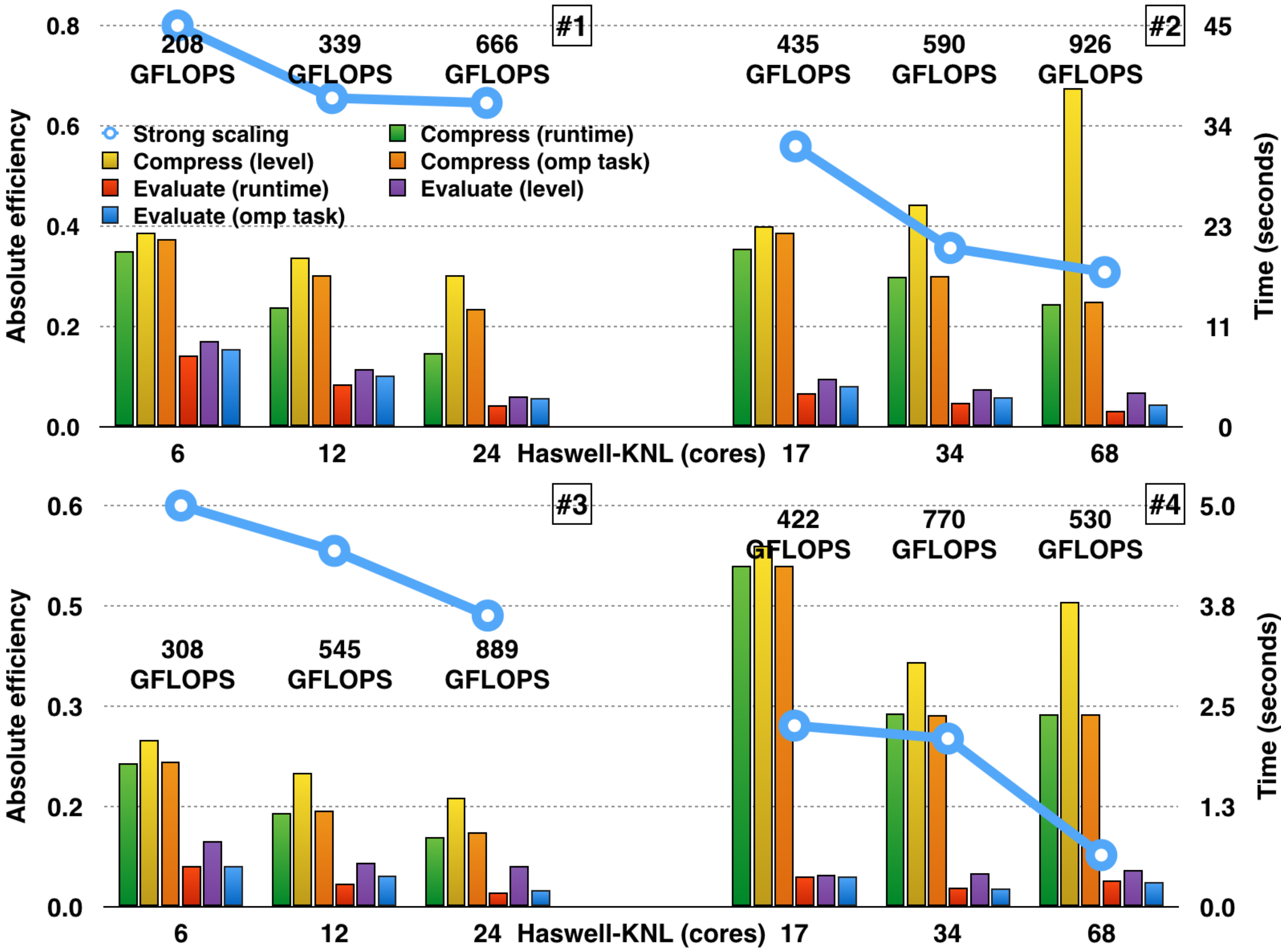}
  \caption{
  Strong scaling on a single Haswell and KNL node
  (y-axis, time in seconds on the right, absolute efficiency to the
  peak \texttt{GFLOPS} on the left).
  We use $s512$, $\tau1E-5$ and $r512$.
  \#\ref{exp:covtypehaswel} and \ref{exp:covtypeknl}
  use \textbf{COVTYPE} to create a Gaussian kernel matrix with 
  $m800$ and $12\%$ budget ($h=0.1$), achieving $\epsilon_2 = \num{2E-3}$
  with average rank $487$.
  \#\ref{exp:k02haswell} and \#\ref{exp:k02knl} use \textbf{K02}
  with $m512$ and $3\%$ budget, achieving $\epsilon_2=\num{5E-5}$
  but only with average rank $35$. 
  We increase the number of cores up to 24 Haswell cores and 68 KNL
  cores. Each set of experiments contains compression time and evaluation time on
  three different parallel schemes: wall-clock time, level-by-level and omp tasks.
  We cannot perform scaling experiments for the hybrid CPU-GPU platform
  (see \tabref{tab:arch} for GPU performance).
  }
  \label{fig:scaling}
\end{figure}
\begin{figure}[!t]
  \centering
  \includegraphics[scale=.24]{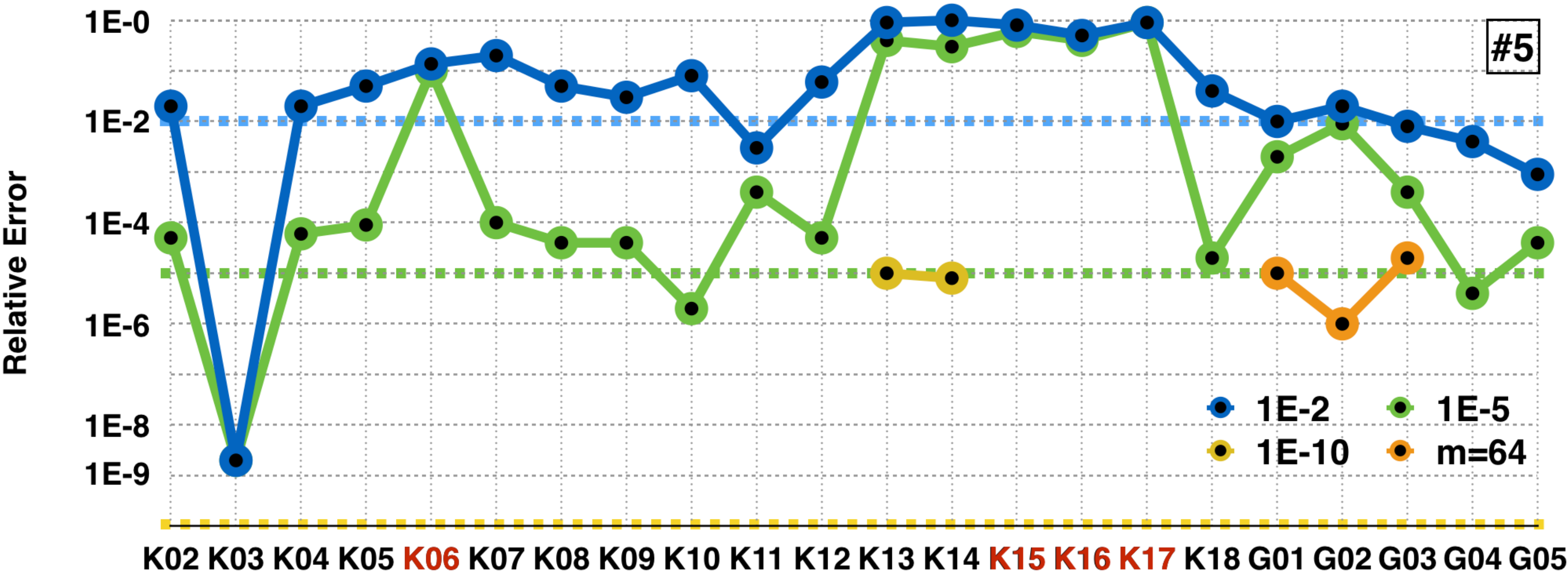}
  \caption{
    \#\ref{exp:robust}, relative error $\epsilon_2$ 
    (y-axis, the smaller the better) 
    on all matrices (x-axis) using angle distance. 
    Blue bars use $\tau\num{1E-2}$ and $1\%$ budget
    (except for \textbf{K6}, \textbf{K15}, \textbf{K16}, \textbf{K17}, other matrices take $0.8$s
    to compress and $0.1$ to evaluate in average).
    Green bars use $\tau\num{1E-5}$ and $3\%$ budget
    (in average, compression takes $1$s and
    evaluation takes $0.2$s).
    Red labels denotes matrices that do not compress.
    \textbf{K13} and \textbf{K14} have hierarchical low-rank structure, but the adaptive ID underestimates the rank.
    \textbf{K13} and \textbf{K14} can reach high accuracy (yellow plots) 
    with $\tau\num{1E-10}$ and $3\%$ budget ($1.0$s in compression
    and $0.2$s in evaluation).
  }
  \label{fig:robust}
\end{figure}
\begin{figure}[!t]
  \centering
  \includegraphics[scale=.29]{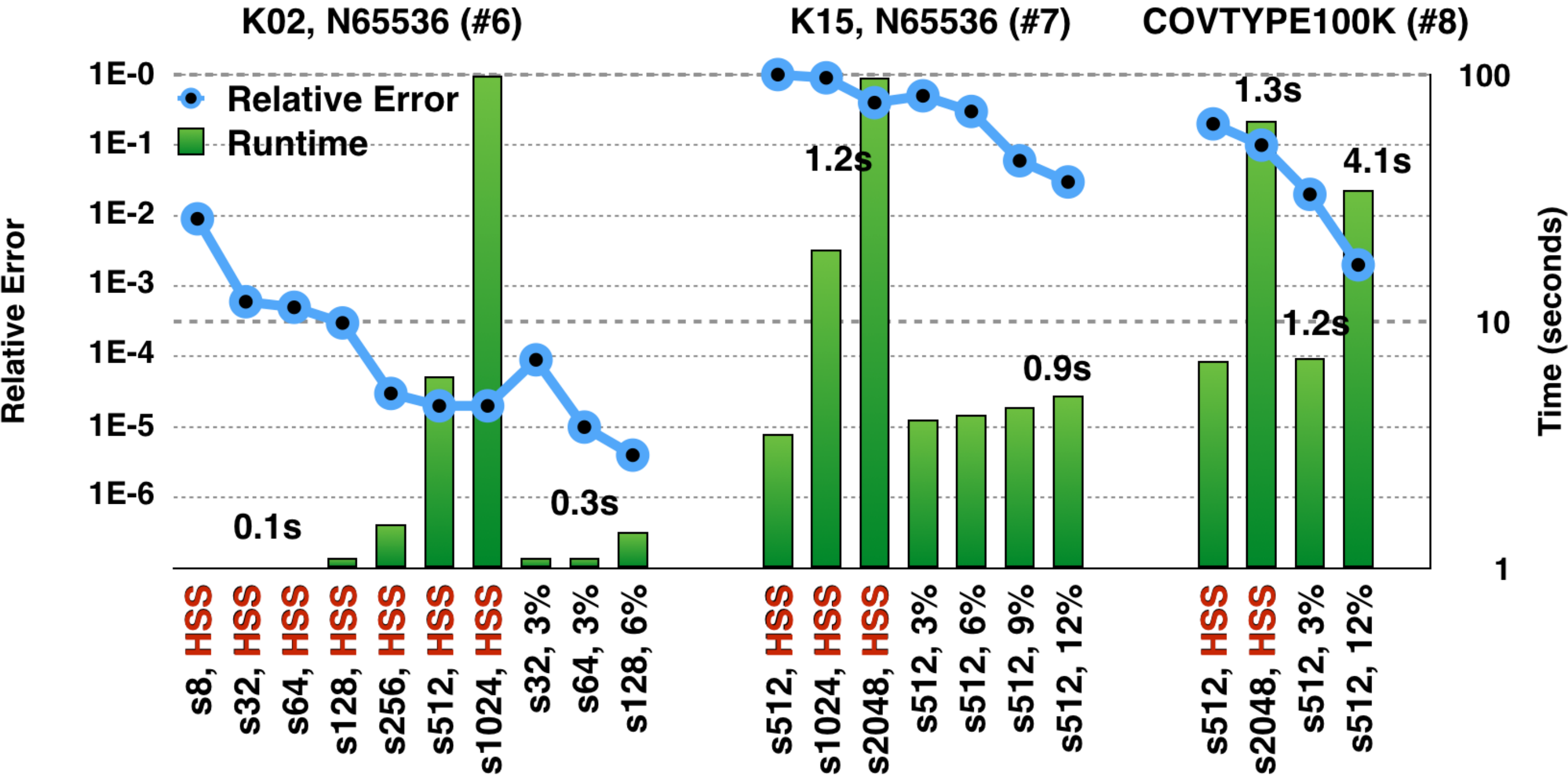}
  \caption{Comparison between HSS and FMM in
  wall-clock time (seconds, green bars, right y-axis) and accuracy ($\epsilon_2$, blue
  plots, left y-axis).
  In \#\ref{exp:k02robust}, \#\ref{exp:k15robust}
  and \#\ref{exp:covtyperobust},
  we use \textbf{K02}, \textbf{K15} ($m512$) and \textbf{COVTYPE} ($m800$)
  datasets. 
  The fixed rank and budget are labeled on x-axis.
  The green bar is the total wall-clock time including compression and evaluation
  on $512$ right hand sides. 
  For some experiments, we also provide wall-clock time for evaluation to contrast
  the trade-off of using high rank and high budget.
  }
  \label{fig:fmmvshss}
\end{figure}
\begin{figure}[!t]
  \centering
  \includegraphics[scale=.23]{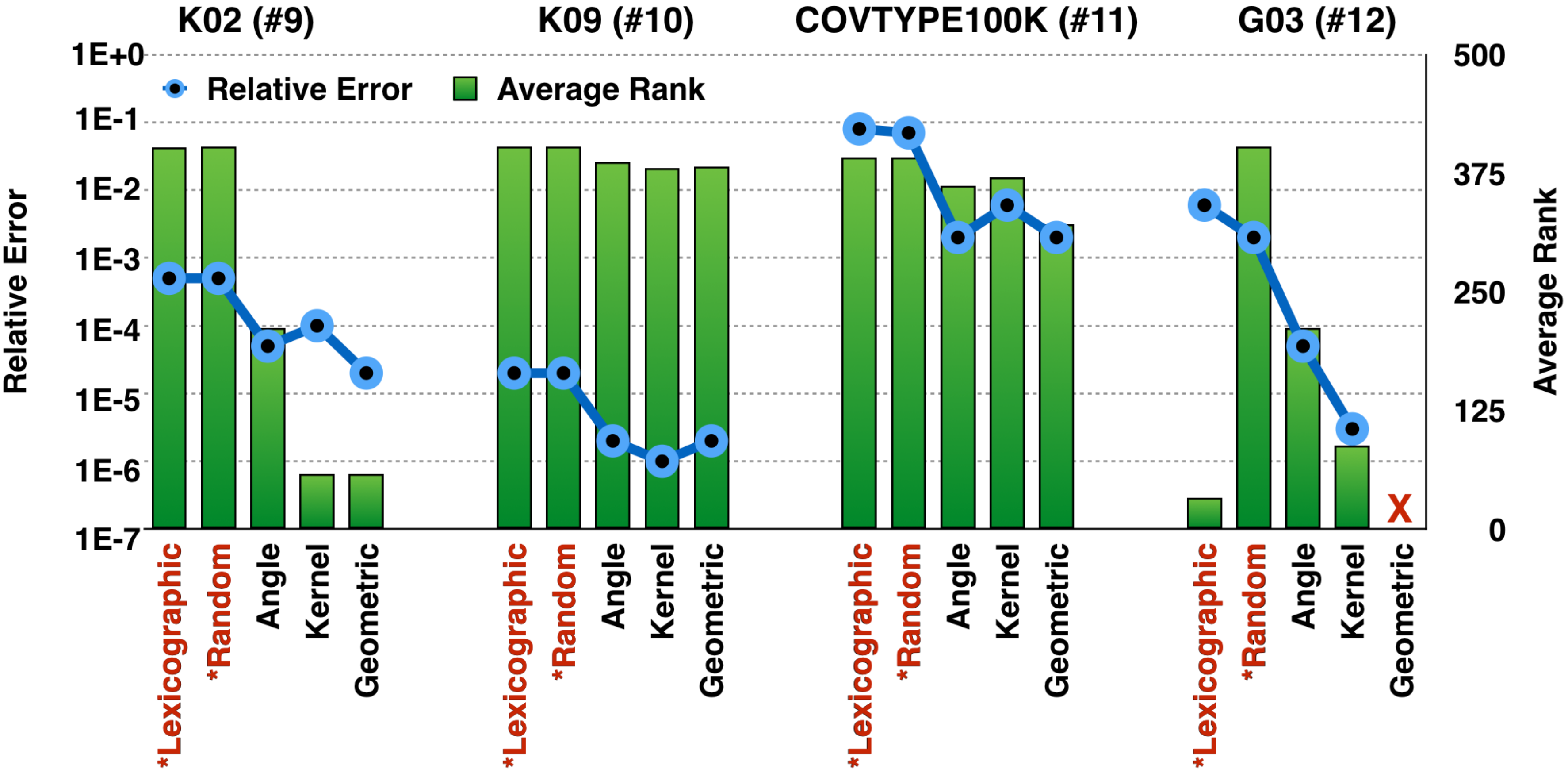}
  \caption{
    Accuracy (left y-axis) and rank (right, x-axis) comparison:
    \textbf{Lexicographic}, \textbf{Random}, \textbf{Kernel 2-norm},
    \textbf{Angle} and \textbf{Geometric}.
    We use $\tau\num{1E-7}$, $s512$, $m64$. For methods that define 
    \emph{distance}, we use $k32$ and $3\%$ budget.
    \textbf{G03} is a graph Laplacian; thus, using
    \textbf{Geometric} distance is impossible.
  }
  \label{fig:partition}
\end{figure}

We label all experiments from \#\ref{exp:covtypehaswel} to
\#\ref{exp:g04knl} in tables and figures.
We perform \ipoint{strong scaling results} on a single Haswell and KNL node
in \figref{fig:scaling}, comparing different scheduling schemes.
In \figref{fig:robust}, we examine the \ipoint{accuracy} of \gofmm{}  for the different matrices; notice that not all 22 matrices admit good hierarchical low-rank
structures in the original order (lexicographic). In \figref{fig:fmmvshss}, we
\ipoint {compare FMM} ($S\neq 0$ in \eqref{e:hmatrix}) \ipoint{to HSS} ($S=0$)
and show an example in which increasing direct evaluations in FMM results in higher accuracy and shorter wall-clock time. 
In \figref{fig:partition}, we present a comparison between 
\ipoint{five permutation} schemes; matrix-defined Gram distances work quite well.

For reference, we compare \gofmm{} to \ipoint{three other codes}:
\texttt{HODLR} and \texttt{STRUMPACK} ($S=0$ in these codes) in
\tabref{tab:existingSoftware} and \texttt{ASKIT} (high-$d$ FMM) in \tabref{tab:askitcomparison}. 
The two first codes do not permute $K$. \texttt{ASKIT} is similar to \texttt{GOFMM} but uses level-by-level traversals, does not produce a symmetric $\tilde{K}$,  and requires points.   Finally, we test \gofmm{} on \ipoint{four different architectures} in \tabref{tab:arch}; the performance of \texttt{GOFMM} correlates with the performance of BLAS/LAPACK.

\textbf{Strong scaling (\figref{fig:scaling}).}
In 
\#\rownumber\label{exp:covtypehaswel}, 
\#\rownumber\label{exp:covtypeknl}, 
\#\rownumber\label{exp:k02haswell}, 
\#\rownumber\label{exp:k02knl}, 
we use a 24-core Haswell and a 68-core KNL to perform strong scaling
experiments. Each set of experiments contains 6 bars including 3 different
parallel schemes on both \algref{a:compress} and \algref{a:evaluate}.
The blue dot indicates the absolute efficiency (ratio to the peak) of 
our \textbf{evaluation} using dynamic scheduling.
\#\ref{exp:covtypehaswel} and
\#\ref{exp:covtypeknl} require $12\%$ budget with 
average rank $487$ to achieve \num{2E-3}. 
This compute-bound problem can reach $65\%$ peak performance on Haswell and $33\%$ on KNL.
However, \#\ref{exp:k02haswell} and \#\ref{exp:k02knl}
only require $3\%$ budget with average rank $35$ 
to achieve $\num{5E-5}$. As a result, this memory-bound problem does not 
scale ($46\%$ and $8\%$\footnote{\scriptsize
The average rank of \#\ref{exp:k02knl} is too small. Except for
\texttt{L2L} tasks, other tasks can only reach about $5\%$
of the peak during the evaluation.
We suspect that MKL' \texttt{SGEMM} uses a $30\times16$ 
micro-kernel to perform a $30\times256\times16$ rank-$k$ update
each time. For an $m\times k\times n$ \texttt{SGEMM} to be efficient, 
$m$ and $n$ usually need to be at least four times of the
micro-kernel size in each way. 
In \#\ref{exp:k02knl}, many \texttt{SGEMM}s have $m<30$. 
Still the micro-kernel must compute $2\times30\times256\times16$ FLOPS.
These sparse FLOPS are not counted in our experiments.
}) very well.  In \#\ref{exp:k02knl}, we can even observe slow down from 34-core to 68-core.
This is because the wall-clock time is bounded by the task in the critical path;
thus, increasing the number of cores does not help.

Throughout, we can observe that the wall-clock time for compression is less than
the level-by-level and \texttt{omp task} traversals. While the work
of \texttt{SKEL} is bounded by $2s^3$, parallel \texttt{GEQP3} in the
level-by-level traversal does not scale (especially on KNL).  On the
other hand, task based implementations can execute \texttt{COEF}
and \texttt{Kba} out-of-order to maintain the parallelism. Our wall-clock time
is better that \texttt{omp task} since we use the cost-estimate model
for scheduling.

\textbf{Accuracy (\figref{fig:robust}).}
We conduct \#\rownumber\label{exp:robust} to examine the accuracy of \gofmm{} (up to single precision).
Given $m512$, $s512$ and $r512$, we report relative error $\epsilon_2$
on {\bf K02-18} and {\bf G01-G05} using the \textbf{Angle} distance with
two tolerances: $\num{1E-2}$ (in blue) and $\num{1E-5}$ (in green). 
Throughout, except for \textbf{K06}, \textbf{K15}--\textbf{K17} (high rank),
\textbf{K13}, \textbf{K14} (underestimating the rank), and
\textbf{G01}--\textbf{G03} (requiring smaller leaf size $m$),
other matrices can usually achieve high accuracy with tolerance $\num{1E-5}$
($0.9$s in compression and $0.2$s in evaluation).
Our adaptive ID underestimates the rank of \textbf{K13} and \textbf{K14} 
such that $\epsilon_2$ is high. By imposing a smaller tolerance 
$\num{1E-10}$ (yellow plots), both matrices 
reach $\num{1E-5}$ ($1$s in compression and $0.2$s in evaluation).
\textbf{K6}, \textbf{K15}--\textbf{K17} have high ranks in the 
off-diagonal blocks; thus they cannot be compressed with $s512$ 
and $3\%$ budget. \textbf{G01}--\textbf{G03} requires direct evaluation in the off-diagonal 
blocks to reach high accuracy. When we reduce the leaf node size from $512$ to $64$, 
we can can still reach $\num{1E-5}$ (orange plots).
However, decreasing leaf size to $64$ results in a longer wall-clock time
($0.8s$ in evaluation), because small $m$ hurts performance.
Overall, we can observe that \texttt{GOFMM} can quite robustly discover
low-rank plus sparse structure from different SPD matrices.
We now investigate how increasing the cost (either with higher rank
or more direct evaluations) can  improve accuracy.

\textbf{Comparison between FMM and HSS (\figref{fig:fmmvshss}).}
We use \#\rownumber\label{exp:k02robust}, \#\rownumber\label{exp:k15robust},
and \#\rownumber\label{exp:covtyperobust} to show that even with more
evaluations, FMM can be faster than HSS for the same accuracy. 
For HSS the relative error in \#\ref{exp:k02robust} (blue plots) plateaus at
\num{5E-4}. Further increasing rank from 256 to 512 (or even 1,024) results
in $O(s^3)$ work (green bars). Using a combination of low-rank
($s64$) and $3\%$ direct evaluation, FMM can achieve higher accuracy
with little increment in the evaluation time (compression time remains
the same).  Similarly, in \#\ref{exp:covtyperobust} we can observe that
by using $s512$ and $3\%$ budget we achieve better accuracy than the
HSS approximation ($s2048$) in less time.

\textbf{Permutations (\figref{fig:partition}).}
Here we test different permutations
(\#\rownumber\label{exp:k02partition},
\#\rownumber\label{exp:k09partition},
\#\rownumber\label{exp:covtypepartition}, and
\#\rownumber\label{exp:g03partition}) 
to discuss the different distances in \gofmm{}.
In each set of experiments, we present relative
error (blue plots) and average rank (green bars)
for five different schemes.
The first two schemes use lexicographic or random order
to recursively permute $K$. 
Since there is no \emph{distance} defined, 
these two schemes can only use HSS approximation.
The \textbf{Angle} and \textbf{Kernel}
\emph{distance} use the corresponding Gram distances~\secref{s:oblivious}.
Finally, we also use standard geometric distance from points. 
For the last three schemes, we use $\kappa32$ and $3\%$
budget. Overall, we can observe that the distance metric
is important in discovering low-rank structure
and improving accuracy.
For example, in \#\ref{exp:k02partition}, \textbf{Kernel} and 
\textbf{Geometric} show much lower average rank
than others.
In \#\ref{exp:k09partition} and \#\ref{exp:covtypepartition},
although the average ranks are not significantly
different, distance-based methods usually have higher accuracy. 
Finally, we observe for matrix \textbf{G03} in \#\ref{exp:g03partition} where no coordinate information exists,
our geometry-oblivious methods can still compress
the matrix. Although the lexicographic permutation has
very low rank, the error is large. This is because
the uniform samples for the low-rank approximation are poor. 
\textbf{Angle} and \textbf{Kernel} distance use neighbors for importantance sampling, 
which greatly improves the quality of the low-rank approximation.

\begin{table}
\input{tables/otherSoftware2.tex}
\caption{
  Wall-clock time comparison (in seconds) between \texttt{HODLR},
  \texttt{STRUMPACK}, and \texttt{GOFMM}.
  For \textbf{K02}--\textbf{K12}, we use $N=36K$. \textbf{K17} uses $N=32K$,
  and  \textbf{G03} uses $N=65K$.
  For all software, we use leaf node size $m512$ and $1024$ right hand sides.
  We control other parameters ($\tau$ and $s$) for each software to target 
  the same relative error (\num{1E-4}).
%
}
\label{tab:existingSoftware}
\end{table}

\begin{table}
\input{tables/askit_comparison.tex}

\caption{
  Wall-clock time (in seconds) and accuracy $\epsilon_2$ comparison with \texttt{ASKIT}.
  For both methods, we use $\kappa=32$, $m = s = 512$ and $r1$.
  \texttt{ASKIT} use the $\tau$ reported in the table, and we adjust
  the tolerance of \texttt{GOFMM} to match the accuracy.
  For all experiments, \texttt{GOFMM} uses $7\%$ budget.
  The amount of direct evaluation performed by \texttt{ASKIT} is decided
  by $\kappa$.
}
\label{tab:askitcomparison}
\end{table}

\textbf{Comparison to existing software (\tabref{tab:existingSoftware}, \tabref{tab:askitcomparison}).}
We compare our methods to 
\texttt{HODLR}~\cite{ambikasaran-darve13}, 
\texttt{STRUMPACK}~\cite{rouet-li-e16}, and 
\texttt{ASKIT}~\cite{march-xiao-yu-biros-sisc16}.
Let us summarize some key differences. \texttt{HODLR} uses the
Adaptive Cross Approximation (ACA, partial pivoted LU) for constructing the low-rank blocks (using the \texttt{Eigen} library). Its evaluation requires $O(N\log{N})$ work since the $U$, $V$ matrices are not nested. 
\texttt{STRUMPACK} constructs an HSS representation
in $O(N\log{N})$ work.
This is done by using a randomized ID according to \cite{liberty2007randomized}. 
We used their black-box compression routine with a uniform random distribution 
and a Householder rank-revealing QR.
Once the matrix is compressed, the evaluation time is $O(N)$
per right hand side. \texttt{STRUMPACK} supports multiple right hand sides.
\texttt{ASKIT}'s FMM evaluation has similar complexity as \texttt{GOFMM},
but the amount of direct evaluation is only decided by $\kappa$.
For \texttt{GOFMM}, we further introduce the budget to restrict the
cost. For all comparisons, we try to match the accuracy by controlling
different parameters ($\tau$, $s$, and $\kappa$). Notice
that \texttt{ASKIT} and \texttt{STRUMPACK} support MPI,
whereas \gofmm{} does not. We have not used MPI for distributed 
environment in our experiments.

In \tabref{tab:existingSoftware}, we target final accuracy $\epsilon_2 =
\num{1E-4}$. 
\texttt{GOFMM} uses \textbf{Angle} distance for neighbor search and
tree partitioning. \texttt{HODLR} and \texttt{STRUMPACK} do not have 
built-in partitioning schemes for dense matrices. 
\texttt{STRUMPACK} fails to compress \textbf{K04} (Gaussian kernel in 6D) and
\textbf{K07} (Laplace kernel in 6D). The first two ones, this is because
the lexicographic order does not admit a good \hmatrix{} approximation. The matrix needs to be permuted.
\textbf{K17} is difficult to compress with a pure hierarchical low-rank
matrix. Finally, \textbf{G03} performs better when $S\neq 0$. 
\texttt{HODLR} and \texttt{STRUMPACK} must increase the off-diagonal 
ranks to match the accuracy and thus the cost increases. \texttt{GOFMM} just uses $S$ and it is about
$25\times$ faster in compression and about $1.5\times$ faster in evaluation.

In \tabref{tab:askitcomparison}, we compare \gofmm{} (with geometric distances) to \texttt{ASKIT}. 
\texttt{ASKIT} uses level-by-level traversals in both compression and evaluation.
Since \texttt{ASKIT} only evaluates a single right hand side,
we use $r=1$. The compression time is inconclusive for
\#\ref{exp:k04askit1}--\#\ref{exp:k04askit4}; the average ranks
used in two methods are quite different. 
The benefit of \emph{out-of-order} traversal appears in
\#\ref{exp:k06askit1}--\#\ref{exp:k06askit4} where both methods reach 
the maximum rank $s$. The speedup in evaluation is not significant. 
\texttt{GOFMM} can get up to $2\times$ speedup in compression. 

\begin{table}[!t]
\setlength\tabcolsep{4.8pt}
\centering
{
  \begin{tabular}{|r|>{\columncolor[gray]{0.8}}r|rr|rr|rr|} 
  \hline
  \rowcolor[gray]{0.8}
  \# & Arch & Budget & $\epsilon_2$ & Comp & GFs & Eval & GFs \\
  \hline
  \rowcolor[gray]{0.8}
  \multicolumn{8}{|c|}{\textbf{MNIST60K}, $h1$, $\kappa32$, $m512$, $s128$, $r256$} \\
  \hline
  \rownumber\label{exp:mnistarm}   & ARM & $5\%$ & 5E-3 & 285 &  3 & 520 & 12 \\
  \hline
  \rowcolor[gray]{0.8}
  \multicolumn{8}{|c|}{\textbf{COVTYPE100K}, $h1$, $\kappa32$, $m512$, $s128$, $r256$} \\
  \hline
  \rownumber\label{exp:covtypearm}   & ARM & $5\%$  & 8E-4 & 71 &  2 &  61 & 10 \\
  \hline
  \rowcolor[gray]{0.8}
  \multicolumn{8}{|c|}{\textbf{COVTYPE100K}, $h0.1$, $\kappa32$,  $m800$, $s512$,
  $r512$} \\
  \hline
  \rownumber\label{exp:covtypecpu} & CPU & $12\%$ & 2E-3 & 30 & 30 & 4.1 &  679 \\
  \rownumber\label{exp:covtypegpu}   & CPU$+$GPU & $12\%$ & 3E-3 & 33 & 29 & 1.7 & 1952  \\
  \rownumber\label{exp:covtypemic}   & KNL & $12\%$ & 2E-3 & 48 & 25 & 3.2 & 1125  \\
  \hline
  \rowcolor[gray]{0.8}
  \multicolumn{8}{|c|}{\textbf{HIGGS500K}, $h0.9$, $\kappa64$, $m1024$, $s256$, $r512$} \\
  \hline
  \rownumber\label{exp:higgscpu} & CPU & $0.3\%$  & 2E-1 & 102 & 18 & 3.3 & 592 \\
  \rownumber\label{exp:higgsgpu}   & CPU$+$GPU & $0.3\%$ & 2E-1 & 180 & 12 & 1.7 & 1147  \\
  \rownumber\label{exp:higgsknl}   & KNL & $0.3\%$ & 2E-1 & 121 & 17 & 2.2 & 872  \\
  \hline
  \rowcolor[gray]{0.8}
  \multicolumn{8}{|c|}{\textbf{K02}, $N65536$, $\kappa32$, $m512$, $s512$, $r512$} \\
  \hline
  \rownumber\label{exp:k02cpu} & CPU & $3\%$  & 9E-5 &      1 & 25 & 0.2 &  889 \\
  \rownumber\label{exp:k02gpu} & CPU$+$GPU & $3\%$ & 1E-4 & 2 & 12 & 0.1 & 2175 \\
  \rownumber\label{exp:k02mic} & KNL & $3\%$ & 1E-4        & 3 & 11 & 0.3 & 530  \\
  \hline
  \rowcolor[gray]{0.8}
  \multicolumn{8}{|c|}{\textbf{K15}, $N65536$, $\kappa32$, $m512$, $s512$, $r1024$} \\
  \hline
  \rownumber\label{exp:k15cpu} & CPU & $10\%$  & 2E-1 &  6.0  & 81 & 1.1  &
  1495  \\
  \rownumber\label{exp:k15gpu} & CPU$+$GPU & $10\%$ & 2E-1  & 7.8  & 62  &
  0.66 &
  2514 \\
  \rownumber\label{exp:k15knl} & KNL & $10\%$ & 2E-1        & 9.2 & 53 & 1.3 & 1549\\
  \hline
  \rowcolor[gray]{0.8}
  \multicolumn{8}{|c|}{\textbf{G03}, $N65536$, $\kappa32$, $m128$, $s512$, $r512$} \\
  \hline
  \rownumber\label{exp:g03cpu} & CPU & $3\%$       & 4E-5 &  4.8  & 37 & 0.5 &
  1122  \\
  \rownumber\label{exp:g03gpu} & CPU$+$GPU & $3\%$ & 3E-5 &  7.9  & 19  & 0.53 & 962 \\
  \rownumber\label{exp:g03knl} & KNL & $3\%$       & 5E-5        & 11.8 & 9.1 &
  0.6 & 741  \\
  \hline
  \rowcolor[gray]{0.8}
  \multicolumn{8}{|c|}{\textbf{G04}, $N89400$, $\kappa32$, $m512$, $s512$, $r512$} \\
  \hline
  \rownumber\label{exp:g04cpu} & CPU & $3\%$       & 4E-6 &  1.8  & 21 & 0.3  &
  787\\
  \rownumber\label{exp:g04gpu} & CPU$+$GPU & $3\%$ & 4E-6 & 4.0  & 10  & 0.13 & 2277 \\
  \rownumber\label{exp:g04knl} & KNL & $3\%$       & 4E-6 & 4.2 & 9 & 1.5 &
  215  \\
  \hline
  \end{tabular}
}
\caption{
  Accuracy $\epsilon_2$, wall-clock time (in seconds) and efficiency (in \texttt{GFLOPS})
  on four architectures. Because our ARM platform only has a 8GB SD card and
  2GB DRAM, we only perform kernel matrices ($K_{ij}$ computed on the fly) 
  with small $r$ and $s$. Note that in the CPU$+$GPU experiment, the compression is run on the CPU (see ~\secref{s:par}).
}
\label{tab:arch}
\end{table}

\textbf{Different architectures.}
In \tabref{tab:arch}, we present wall-clock time and \texttt{GFLOPS}
of \gofmm{} on four architectures for different problems.
We want to show that the efficiency of \texttt{GOFMM} is
portable and only relies on BLAS/LAPACK libraries.  

In \#\ref{exp:mnistarm} and \#\ref{exp:covtypearm},
we show that a quad-core ARM processor can handle
up to 100K fast matrix-multiplication. Because we only have limited 
memory (2GB) and storage (8GB), in \texttt{GOFMM} we
compute $K_{ij}$ on the fly (in detail, we compute
$K_{\beta\alpha}$ with a \texttt{GEMM} using the 2-norm expansion).
\#\ref{exp:mnistarm} takes much longer than \#\ref{exp:covtypearm}
because the cost of evaluating $K_{ij}$ is proportional to the point dimensions
of the dataset (\textbf{MNIST} in 780D and \texttt{COVTYPE} in 54D).
Because there is no active cooling, the ARM gets overheated and is
forced to reduce its clockrate. That is why we can only reach $30\%$ of peak
during the evaluation. 

Experiments \#\ref{exp:covtypecpu} to \#\ref{exp:higgsknl} are computed
in double precision. 
With $12\%$ budget,
our evaluation can reach $68\%$ peak performance on Haswell, $37\%$ on KNL and $38\%$ 
on a hybrid Haswell-P100 system.
The performance degrades in \#\ref{exp:higgscpu}--\ref{exp:higgsknl} because
the rank is limited to 256, and $0.3\%$ direct evaluation 
is not enough to create large \texttt{GEMM} calls.
For kernel matrices, the \texttt{GFLOPS} for compression are usually
higher because computing $K_{ij}$ requires floating point operations.
For example, compression of \textbf{COVTYPE} (in 54D) has higher 
\texttt{GFLOPS} than \textbf{HIGGS} (in 28D).
This is not only because \textbf{COVTYPE} is a dataset with high 
dimensionality, but we also use a higher rank $s512$ such that 
\texttt{GEQP3} and \texttt{TRSM} can be more efficient.

Finally, we present performance results on several matrices 
(\#\ref{exp:k02cpu}--\ref{exp:g04knl}) in single precision.
With $10\%$ budget in \textbf{K15},
our evaluation can reach $75\%$ peak on Haswell, $25\%$ on KNL and $25\%$ 
on a hybrid Haswell-P100 system.
This performance requires large leaf node size $m$ and sufficient
direct evaluations (e.g. \#\ref{exp:k02cpu}--\#\ref{exp:g04knl}).
Since \textbf{G03} requires small $m$, our \texttt{GFLOPS} efficiency
degrades due to the dependency on the BLAS/LAPACK routines. 
Notice that $m128$ is not large enough for \texttt{GEMM}
to reach high performance on KNL and GPUs.
For \textbf{G04}, we use $m512$ but KNL (\#\ref{exp:g04knl}) does 
not perform very well. 
The same problem occurs in \figref{fig:scaling}: 
the average rank in \textbf{G04} is too small. 
Additionally, we do not observe huge performance degradation on GPUs
(\#\ref{exp:g04gpu}). 
This is because we enforce our scheduler to schedule \texttt{L2L} tasks
to the GPU; thus, tasks with small ranks (\texttt{N2S} and
\texttt{S2N}) are mostly consumed by the host CPU.
The comparison between \#\ref{exp:g04gpu} and \#\ref{exp:g04knl} 
is a good example that highlights the goal of heterogeneous parallel
architectures. CPUs with short vector lengths are suitable for 
tasks with very low ranks  (\texttt{N2S} and \texttt{S2N}). 
On the contrary, GPUs are the method of choice for \texttt{FLOPS} intensive tasks (\texttt{L2L}).
We cannot solve such problems with only one architecture  efficiently.

%% file: tables/otherSoftware2.tex
\centering \small 
\setlength{\tabcolsep}{3pt}
\begin{tabular}{|r|>{\columncolor[gray]{0.8}}L|rrr|rrr|rrr|} 
\hline 
&& \multicolumn{3}{c|}{\texttt{HODLR}} & 
   \multicolumn{3}{c|}{\texttt{STRUMPACK}} &
   \multicolumn{3}{c|}{\texttt{GOFMM}} \\
\hline
\rowcolor[gray]{0.8}
\# & case  & 
$\epsilon_2$ & {Comp} & {Eval} & 
$\epsilon_2$ & {Comp} & {Eval} & 
$\epsilon_2$ & {Comp} & {Eval} \\ 
\hline
\rownumber\label{exp:k02cmp} & K02 & 
\accnum{6.26302e-05} & \num{0.62} & \num{2.74} & 
\accnum{0.000112964} & \num{9.177} & \num{0.559} & 
\accnum{1.523281E-05} & \num{1.0} & \num{0.275162} \\
\rownumber\label{exp:k04cmp} & K04 & \accnum{6.40027e-05} & \num{0.67} & \num{2.69} & \accnum{0.000120927} & \num{507.754} & \num{7.8344} & \accnum{1.810217E-05} & \num{0.994} & \num{0.51}\\
\rownumber\label{exp:k07cmp} & K07 & \accnum{7.26224e-05} & \num{0.87} & \num{3.06} & \accnum{0.000227985} & \num{528.377} & \num{8.20218} & \accnum{3.933854E-05} & \num{0.588} & \num{0.240} \\
\rownumber\label{exp:k12cmp} & K12 & 
\accnum{6.13306e-05} & \num{0.67} & \num{2.65} & 
\accnum{0.000196841} & \num{18.79} & \num{0.75} & 
\accnum{1.159227E-04} & \num{0.6} & \num{0.220}  \\
\rownumber\label{exp:k17cmp} & K17 & 
\accnum{0.100087} & \num{862.16} & \num{37.59} & 
\accnum{0.172307} & \num{663.399} & \num{8.17784} & 
\accnum{9.440622E-02} & \num{48.8} & \num{3.139} \\
\rownumber\label{exp:g03cmp} & G03 & 
\accnum{0.000270227} & \num{12.92} & \num{9.67} & 
\accnum{0.0254869} & \num{29.7846} & \num{1.3314} & 
\accnum{8.233808E-05} & \num{0.5} & \num{0.786} \\
\hline 
\end{tabular}

%% file: tables/askit_comparison.tex
\centering \small 
\setlength{\tabcolsep}{3pt}
\begin{tabular}{|r|>{\columncolor[gray]{0.8}}r|rr|rrr|rrr|} 
\hline 
& \multicolumn{3}{c|}{Parameters} & 
   \multicolumn{3}{c|}{\texttt{ASKIT}} &
   \multicolumn{3}{c|}{\texttt{GOFMM}} \\
\hline
\rowcolor[gray]{0.8}
\# & 
case & $N$ & $\tau$ & 
$\epsilon_2$ & {Comp} & {Eval} & 
$\epsilon_2$ & {Comp} & {Eval} \\ 
\hline
\rownumber\label{exp:k04askit1} &
  K04 & \num{36864} & \accnum{1e-03} & 
  \accnum{0.000184409} & \num{0.2878} & \accnum{0.021256} &
  \accnum{2.362880E-04} & \num{0.5908} & \accnum{1.782e-02} \\
\rownumber\label{exp:k04askit2} &
  K04 & \num{36864} & \accnum{1e-06} & 
  \accnum{7.51693e-07} & \num{1.4374392} & \accnum{0.0388908} &
  \accnum{6.715137E-07} & \num{0.9557} & \accnum{2.516e-02} \\
\rownumber\label{exp:k04askit3} &
  K04 & \num{65536} & \accnum{1e-03} & 
  \accnum{0.00020665} & \num{1.013103} & \accnum{0.043227} &
  \accnum{1.804782E-04} & \num{1.1819} & \accnum{4.217e-02} \\
\rownumber\label{exp:k04askit4} &
  K04 & \num{65536} & \accnum{1e-06} & 
  \accnum{7.35954e-07} & \num{2.2384459} & \accnum{0.0834641} &
  \accnum{6.132387E-07} & \num{1.6681} & \accnum{4.403e-02} \\
\hline
\rownumber\label{exp:k06askit1} &
  K06 & \num{36864} & \accnum{1e-03} & 
  \accnum{0.0351756} & \num{6.556525} & \accnum{0.058805} &
  \accnum{3.428371E-02} & \num{3.3116} & \accnum{3.939e-02} \\
\rownumber\label{exp:k06askit2} &
  K06 & \num{36864} & \accnum{1e-06} & 
  \accnum{0.0211519} & \num{7.4169631} & \accnum{0.0563569} &
  \accnum{3.155044E-02} & \num{4.7576} & \accnum{4.683e-02} \\
\rownumber\label{exp:k06askit3} &
  K06 & \num{65536} & \accnum{1e-03} & 
  \accnum{0.0380554} & \num{11.056666} & \accnum{0.109234} &
  \accnum{4.461605E-02} & \num{5.6683} & \accnum{8.119e-02} \\
\rownumber\label{exp:k06askit4} &
  K06 & \num{65536} & \accnum{1e-06} & 
  \accnum{0.0499491} & \num{11.972589} & \accnum{0.119011} &
  \accnum{4.116939E-02} & \num{7.7440} & \accnum{8.554e-02} \\
\hline 
\end{tabular}

%% file: conclusion.tex
By using the Gramian vector space for SPD matrices, we defined distances between rows of $K$ using only  matrix
values. Using the distances, we introduced \gofmm{} and \hmatrix{} scheme that
can be used to compress arbitrary  SPD matrices (but without accuracy
    guarantees).  In \gofmm{} we use a shared-memory runtime system that
performs \emph{out-of-order} scheduling in parallel to resolve the dynamic
workload due to  
adaptive ranks and the parallelism-diminishing issue during tree traversals.
Our future work will focus on the distributed algorithms and the hierarchical matrix factorization based on our method. 
We also plan to improve the sampling and pruning quality and to reduce the
number of parameters that users need to provide.

%% file: appendix.tex
\subsection{Abstract}
This artifact description appendix comprises the source code, datasets and installation
instruction on a GitHub repository that will be open source and
used to reproduce results for our SC'17 paper.
We also provide all hardware and software configuration in 
\secref{s:sup_setup}.
Due to the double-blind peer reviewing policy, we can only provide 
the url to the repository upon the acceptance.

\subsection{Description} \label{s:sup_setup}

\textbf{Check-list.} 
We briefly describe all meta information.
This program implements an algebraic Fast Multipole Method
with geometric oblivious technique that generalizes to 
SPD matrices.
\begin{itemize}[leftmargin=*]\zapspace
  \item \textbf{Program.} \texttt{GOFMM} is developed in \texttt{C++} (with
    \texttt{C++11} features)
  and \texttt{CUDA}, employing \texttt{OpenMP} 
  for shared memory parallelism using a self-contained runtime system.
\item \textbf{Hardware.}We conducted experiments on
        Lonestar5 (two 12-core, 2.6GHz, Xeon E5-2690 v3  ``Haswell''
        per node) and Stampede (68-core, 1.4GHz, Xeon Phi 7250 ``KNL'' per node)
        clusters at the Texas Advanced Computing Center,
Piz Daint (12-core, 2.3GHz, Xeon  E5-2650 v3 and NVIDIA Tesla P100)
at Swiss National Supercomputing Centre,
and an Intrinsyc Open-Q 820 Development Kit (quad-core, 2.2GHz Qualcomm Kyro).
\item \textbf{Compilation.} All software (including \texttt{HODLR}, \texttt{STRUMPACK} and \texttt{ASKIT}) 
are compiled with 
\texttt{intel-16.0 -O3}  on
Lonestar5 and Piz Daint. Stampede uses \texttt{intel-17.0 -O3 -xMIC-AVX512}.
The GPU part uses \texttt{nvcc-8.0 -O3 -arch$=$sm\_60}. 
For Open-Q 820, we cross compile our software with NDK using 
\texttt{gcc-4.9 -O3}.
All CPU and KNL BLAS/LAPACK routines use MKL.
GPU BLAS routines use CUBLAS; on ARM we use QSML
(Qualcomm Snapdragon Math Library).
KNL experiments use Cache-Quadrant configuration. 
\texttt{OpenMP} uses $\texttt{OMP\_PROC\_BIND=spread}$.
\item \textbf{Datasets.} Our 22 matrices \textbf{K02}--\textbf{G05} can be generated using 
  MATLAB scripts (provided in the repo). 
  The urls of the five graphs \textbf{G01}--\textbf{G05} and the
  For real world datasets, we provide urls in \secref{s:setup}.
\item \textbf{Output.} Runtime and total \texttt{FLOPS} of the compression
  and evaluation phase, accuracy $\epsilon_2$ of the first 10 entries
  and the average of 100 entries.
\item \textbf{Experiment workflow.} \texttt{git clone} projects;
  generate datasets; run test scripts; observe the results;
\end{itemize}

\textbf{How delivered.}
Upon acceptance, we will provide the url to the \gofmm{} 
repository on GitHub. The software comprises code, build, and evaluation
instructions, and is provided under GPL-3.0 license.


\textbf{Hardware dependencies.}
For adequate reproducibility, we suggest that reproducers use
the same environment as mentioned above. Notice that we report
absolute \texttt{GFLOPS} and the ratio to the peak performance
in the paper. For approximately reproducing the same results
on a different environment, reproducer should look for platform
that has similar capability.
The theoretical peak performance\footnote{\scriptsize{We estimate the peak 
according to the clockrate
and the \texttt{FMA} throughput. For 24 Haswell cores,
$998=2\times12\times2.6\times16$. For 68 KNL cores,
$3046=68\times1.4\times32$. For 4 ARM cores, $35.2=4\times 2.2\times 4$.
The peak of P100 is reported as 4.7 \texttt{TFLOPS}.
As a reference, MKL \texttt{GEMM} can achieve $87\%$ on a Haswell node and 
$69\%$ on a KNL node. QSML \texttt{GEMM} can achieve $89\%$ on Open-Q 820.
\texttt{cublasXgemm} can achieve $95\%$ on P100.
We assume two KNL VPUs can dual issue \texttt{DFMA}s~\cite{sodani2016knights}.
However, Intel processors may have a different frequency while fully issuing 
  \texttt{FMA}, and the clockrate may drop to 1.0 GHz. This may be the reason
why MKL \texttt{DGEMM} can only achieve 2.1 TFLOPS on KNL.}} in double precision is $998$ GFLOPS per Haswell node, $3,046$ GFLOPS
per KNL node, ($4,700+416$) GFLOPS per Tesla P100 node, and $35.2$ GFLOPS 
per Open-Q820.
The peak \texttt{GFLOPS} doubles for single precision computations.

\textbf{Software dependencies.}
Compilation requires \texttt{C/C++} compilers that support \texttt{c++11}
features and \texttt{OpenMP}. \gofmm{} also requires full functionality 
of BLAS and LAPACK routines.

\subsection{Installation} \label{s:sup_install}
Given the repository url, you should be able to clone the \texttt{master}
branch of the repository. The first step is to edit \texttt{set\_env.sh}
to select the proper compiler and architecture.

\
\begin{verbnobox}[\small]
export GOFMM_USE_INTEL = true  
export GOFMM_USE_INTEL = false 
export GOFMM_USE_CUDA  = true  
export GOFMM_MIC_AVX512= true  
\end{verbnobox}
\

If user want to compile the CUDA code for the hybrid CPU-GPU implementation
then the following variables have to be exported.

\
\begin{verbnobox}[\small]
export GOFMM_GPU_ARCH_MAJOR=gpu
export GOFMM_GPU_ARCH_MINOR=pascal
\end{verbnobox}
\

There are three options for the host (ARM, x86-64 or KNL). Users must 
choose at least one major and minor architecture to compile.
This can be arm/armv8a, x86-64/haswell or mic/knl.

\
\begin{verbnobox}[\small]
export GOFMM_ARCH_MAJOR=arm
export GOFMM_ARCH_MINOR=armv8a

export GOFMM_ARCH_MAJOR=x86_64
export GOFMM_ARCH_MINOR=haswell

export GOFMM_ARCH_MAJOR=mic
export GOFMM_ARCH_MINOR=knl
\end{verbnobox}
\

Although we use \texttt{cmake} to identify BLAS/LAPACK libraies, but we
suggest that user manually setup the path using

\
\begin{verbnobox}[\small]
export GOFMM_QSML_DIR = /path-to-qsml
export GOFMM_MKL_DIR = /path-to-mkl
export GOFMM_CUDA_DIR  = /path-to-cudatoolkit
\end{verbnobox}
\

Finally, users must setup the OpenMP option to enable parallel implementation.
Here for example, we use 68 threads for KNL and \texttt{spread}
\texttt{OpenMP} thread binding. 

\
\begin{verbnobox}[\small]
export OMP_PROC_BIND=spread
export OMP_NUM_THREADS=68
\end{verbnobox}
\

With all these options setup, now we use \texttt{cmake} for compilation. Users can use the following commends.

\
\begin{verbnobox}[\small]
source set_env.sh
mkdir build
cd build
cmake ..
make
make install
cd bin
./run_gofmm_x86
./run_gofmm_gpu
./run_gofmm_knl
\end{verbnobox}
\

\textbf{Cross compilation.}
If your ARM runs with OS that has native compiler and \texttt{cmake} support, 
then the installation instructions above should work just fine. However, while your
target runs an Android OS, which currently does not have a native
\texttt{C/C++} compiler, you will need to cross compile this software on your 
Linux or OSX first.
Although there are many ways to do cross compilation, we suggest that users
follow these instructions:
\begin{itemize}[leftmargin=*]\zapspace
  \item Install Android Studio with LLDB, cmake and NDK support.
  \item Create stand-alone-toolchain from NDK.
  \item Install adb (Android Debug Bridge)
  \item Compile with cmake. It will look for your arm \texttt{gcc/g++},
    \texttt{ar} and \texttt{ranlib} support.
  \item Use the following instructions to push executables and scripts in
    /build/bin to the Android ARM device.
\end{itemize}

\
\begin{verbnobox}[\small]
adb devices
adb push /build/bin/* /data/local/tmp
adb shell
cd /data/local/tmp
./run_gofmm_arm.sh
\end{verbnobox}
\

\subsection{Dataset}
The 22 matrices we use can be generated using MATLAB scripts with the 
corresponding coordinates or graphs. The five graphs we use in the paper
are reported in \secref{s:setup}, and the MATLAB scripts will be provided
as parts of the source code upon acceptance.

\subsection{Experiment workflow}
With the repository url, \texttt{git clone} projects.
Generate datasets using the provided MATLAB script.
Compile \gofmm{} with the instructions in \secref{s:sup_install}.
Run the test script for each architecture. Observe the results.

%
%
%
%
%
%

\subsection{Evaluation and expected result}
For x86-64, ARM and KNL execution, the program will start from
the iterative ANN. The accuracy is reported in every iteration.
Once the neighbor search is done (or skipped), the metric ball
tree partitioning follows.
The program reports runtime and total \texttt{FLOPS} of the compression
  and evaluation phase.
Finally, the accuracy $\epsilon_2$ is reported in two parts: 
the error of the first 10 entries, and the average error of 100 entries.
Notice that in a CPU-GPU hybrid environment, \gofmm{} will
first try to detect the available GPU device. 
If successful, the device name and the available global memory
size should be displayed.
The rest of the execution is the same as our architectures. 